\documentclass[12pt]{article}
\usepackage{amssymb}
\usepackage{amsmath}
\usepackage{amsthm}
\usepackage{amsopn}
\newtheorem{theorem}{Theorem}
\newtheorem{proposition}{Proposition}[section]
\newtheorem{corollary}[proposition]{Corollary}
\newtheorem{lemma}[proposition]{Lemma}

\newtheorem{conjecture}[proposition]{Conjecture}

\let\ds=\displaystyle
\def\ant{S}
\def\coun{\varepsilon}
\def\bn{\begin{equation}}
\def\ed{\end{equation}}

\newcommand{\ein}[2]{e_{#1}[#2]}
\newcommand{\fin}[2]{f_{#1}[#2]}
\newcommand{\hin}[2]{h_{#1}[#2]}

\newcommand{\han}[2]{h_{#1}[#2]}

\def\xp{e}
\def\xm{f}
\def\<{\langle}
\def\>{\rangle}

\newcommand{\CC}{{\mathbb C}}
\newcommand{\NN}{{\mathbb N}}

\newcommand{\ZZ}{{\mathbb Z}}

\let\z=\ZZ
\def\r#1{(\ref{#1})}
\let\rf=\r
\def\ot{\otimes}
\def\a{\alpha}  \def\d{\delta}  
\def\b{\beta} \def\ga{\gamma}

\def\al{\alpha}
\def\Uqg{{U_{q}^{}(\widehat{\mathfrak{g}})}}
\def\Uqgp{{U_{q}^{\prime}(\widehat{\mathfrak{g}})}}
\let\Uqtri=\Uqg
\newcommand{\Uqdva}{{U_q(\widehat{\mathfrak{sl}}_{2})}}
\newcommand{\Uqdvaprim}{{U'_q(\widehat{\mathfrak{sl}}_{2})}}
\def\sk#1{\left(#1\right)}

\def\i{\iota}
\def\on{\operatorname}
\def\l{\lambda}
\def\ggg{\mathfrak{g}}

\def\bb{{\mathfrak b}}
\def\nn{{\mathfrak n}}

\def\C{C}
\def\Uqgln{U_q(\widehat{\mathfrak{gl}}_N)}
\def\Uqgl#1{U_q(\widehat{\mathfrak{gl}}_{#1})}

\def\ord{\prec}
\def\gfr#1#2{\genfrac{}{}{0pt}{}{#1}{#2}}

\def\Gam{\Sigma}
\def\Gama{\Sigma_+}

\def\h{\mathfrak{h}}
\begin{document}
\begin{center}

\bigskip\bigskip

\hfill ITEP-TH-40/05\\
\hfill math/0610398\\
\bigskip
{\Large\bf Weight functions and Drinfeld currents}
\end{center}
\bigskip
\begin{center}
{\bf
B. Enriquez$^{*}$\footnote{E-mail: enriquez@math.u-strasbg.fr},\ \
S. Khoroshkin$^{\star}$\footnote{E-mail: khor@itep.ru},\ \
S. Pakuliak$^{\star\bullet}$\footnote{E-mail: pakuliak@theor.jinr.ru}}\\
\bigskip
$^{*}$
{\it IRMA (CNRS), 7 rue Ren\'e Descartes, F-67084 Strasbourg, France}\\
$^\star$
{\it Institute of Theoretical \& Experimental Physics, 117259
Moscow, Russia}\\
$^\bullet$
{\it Laboratory of Theoretical Physics, JINR,
141980 Dubna, Moscow reg., Russia}\\
\bigskip
\bigskip

{Revised \today}
\end{center}
\begin{abstract}
A universal weight function for a quantum affine algebra is a family of
functions with values in a quotient of its Borel subalgebra,
satisfying certain coalgebraic properties. In representations of the
quantum affine algebra it gives off-shell Bethe vectors and is used
in the construction of solutions of the qKZ equations. We construct a
universal weight function for each untwisted quantum
affine algebra, using projections onto the intersection of Borel subalgebras
of different types, and study its functional properties.
\end{abstract}


\section{Introduction}
The first step of the nested Bethe ansatz method (\cite{KR83}) consists
in the construction of certain rational functions with values in a
representation of a quantum affine algebra or its rational or elliptic
analogue. In the case of the quantum affine algebra
$U_q(\widehat{\mathfrak{gl}}_2)$ these rational functions,
known as the off-shell Bethe vectors, have the form $B(z_1)\cdots B(z_n)v$,
where $B(u)=T_{12}(u)$ is an element of the monodromy matrix (this is a
generating series for elements in the algebra) and $v$ is a highest weight
vector of a finite dimensional representation of $\Uqgl2$.
For the quantum affine algebra $U_q(\widehat{\mathfrak{gl}}_N)$, the off-shell
Bethe vectors are constructed in \cite{KR83} by an inductive procedure
(the induction is over $N$). These Bethe vectors were then used (under the name
`weight functions') in the construction of solutions of the $q$-difference
Knizhnik-Zamolodchikov equation (\cite{VT,S}).
Inductive procedures for the construction of
Bethe vectors were also used in rational models (where the underlying
symmetry algebra is a Yangian) in \cite{BF} and \cite{ABFR}, for
${\mathfrak g} = \mathfrak{sl}_N$; in these cases the Bethe vectors were
expressed explicitly in the quasi-classical limit or using the Drinfeld
twist (see Section \ref{sect:rel}). Bethe vectors for rational models
with ${\mathfrak g} = \mathfrak{o}(n)$ and $\mathfrak{sp}(2k)$ were studied
in \cite{R85}; the twisted affine case $A_2^{(2)}$ was treated in
\cite{T}.

Despite their complicated inductive definition, the weight functions
enjoy nice properties, which do not depend on induction steps. These
are coalgebraic properties, which relate the weight function in a
tensor product of representations with weight functions in the tensor
components (\cite{VT}).

The goal of this paper is to give a direct construction of weight functions,
independent of inductive procedures. For this purpose we introduce the
notion of a universal weight function. This is a family of formal Laurent
series with values in a quotient of the Borel subalgebra of the quantum
affine algebra, satisfying certain coalgebraic properties.
The action of the universal weight function on a highest weight vector
defines a weight function with values in a representation of the
quantum affine algebra, which enjoys coalgebraic properties as in \cite{VT}.

It is well-known that quantum affine algebras, as well as affine Kac-Moody
Lie algebras, admit two different realizations (\cite{D88}). In the first
realization, the quantum affine algebra is generated by Chevalley generators,
satisfying $q$-analogues of the defining relations for Kac-Moody Lie algebras
(\cite{D}). In the second realization (\cite{D88}), generators are the
components of the Drinfeld currents, and the relations are deformations
of the loop algebra presentation of the affine Lie algebra. The quantum
affine algebra is equipped with two coproducts (`standard' and `Drinfeld'),
each of which expresses
simply in the corresponding realization. Both realizations are related to
weight functions: on the one hand, the weight function satisfies coalgebraic
properties with respect to the standard coproduct structure; on the other hand, the notion of a highest weight vector
is understood in the sense of the `Drinfeld currents' presentation.

Our construction of a universal weight function is based on the use of deep
relations between the two realizations. This connection was done in several
steps. The isomorphism between the algebra structures of both sides was
proved in \cite{D88} and \cite{DF}. The coalgebra structures
were then related in \cite{KhT2,KhT3,Be,Dam,DKhP1}; there it was
proved that the `standard' and `Drinfeld' coproducts are related by a twist,
which occurs as a factor in the decomposition of the $R$-matrix for the
`standard' coproduct. A further description of this twist
(in the spirit of the Riemann problem in complex analysis) was
suggested in \cite{ER} and developed in \cite{EF,E2,DKhP}.
Each realization determines a decomposition of the algebra as the product
of two opposite Borel subalgebras (so there are 4 Borel subalgebras).
In its turn, each Borel subalgebra decomposes as the product of
its intersections with the two Borel subalgebras of the other type,
and determines two projection operators which map it to these intersections.
The twist is then equal to the image of
the tensor of the bialgebra pairing between opposite Borel
subalgebras by the tensor product of opposite projections.

In this paper, we give a new proof of these results. For this, we prove a
general result on twists of the double of a finite dimensional Hopf algebra
$A$, arising from the decomposition of $A$ as a product of coideals
(Subection 2.1); this result has a graded analogue (Subsection 2.2),
which can be applied e.g. to quantum Kac-Moody algebras
with their standard coproduct. More importantly, this result has a
topological version (Subsection 2.3). In Section 3,
we show how this topological version implies that the Drinfeld and
standard coproducts are related by the announced twist.

The main result of this paper is Theorem \ref{th1}. It says that the
collection of images of products of Drinfeld currents by the projection
defines a universal weight function. This allows to compute the
weight function explicitly when $\mathfrak{g} = \mathfrak{sl}_2$ or
$\mathfrak{sl}_3$ (see \cite{KhP1}), or when $q=1$ (see identity
(\ref{classical:W})), using techniques of complex analysis and conformal
algebras (\cite{DKh}). We give a conjecture on the general form of the
universal weight function. We describe functional properties of the
weight functions at the formal level, using techniques of (\cite{E}).
We also prove more precise rationality results in two cases:
(a) in the case of finite dimensional modules, as a
consequence of the conjecture on the form of the universal weight
function; (b) in the case of lowest weight modules, the rationality
follows (unconditionally) from a grading argument.

The paper is organized as follows. In Section 2, we prove results on
twists of doubles of Hopf algebras. In Section 3, we recall the definition
of the untwisted quantum affine algebra $\Uqg$, of its coproducts, the
construction of the Cartan-Weyl basis and its relation with the currents
realization of $\Uqg$, following \cite{KhT2}; we also reprove the twist
relation between the two (Drinfeld and standard) coproducts, using Section 2.
In particular, we introduce Borel subalgebras of different types and the
related projection operators. Their definition, relies on a generalization
of the convexity property of the Cartan-Weyl generators to `circular'
Cartan-Weyl generators, see \cite{KhT}; their properties are proved in the
Appendix. their In Section 4, we define and construct universal
weight functions, and prove the main theorem. As a corollary, we derive
analytical properties of our weight functions. In Section 5 we identify
them, in the case of $\Uqdva$, with the expressions familiar in the algebraic
Bethe ansatz theory.

\section{Twists of doubles of Hopf algebras}
\label{sec2.3}
\setcounter{footnote}{0}

\subsection{The finite dimensional case} \label{sect:fd}

Let $A$ be a finite dimensional Hopf algebra.
Assume that $A_1,A_2$ are
subalgebras of $A$
such that:\footnote{If $X$ is a
Hopf algebra, we denote by
$m_X,\Delta_X,S_X,\varepsilon_X,1_X$ its operations.}
(a) the map $m_A : A_1 \otimes A_2 \to A$ is a vector space isomorphism,
(b) $A_1$ (resp., $A_2$) is a left (resp., right)
coideal of $A$, i.e., $\Delta_A(A_1) \subset A \otimes A_1, \quad
\Delta_A(A_2) \subset A_2 \otimes A$.

Let $P_i : A \to A_i$ be the linear maps such that
$P_1(a_1a_2) = a_1 \varepsilon_A(a_2)$, $P_2(a_1a_2)
= \varepsilon_A(a_1)a_2$ for $a_i\in A_i$.
Then we have $m_A \circ (P_1\otimes P_2) \circ \Delta_A = \on{id}_A$.

Let $D$ be the double of $A$
and let $R \in D^{\otimes 2}$ be its
$R$-matrix. Set $R_i := (P_i \otimes \on{id})(R)$. The above identity,
together with $(\Delta_A \otimes \on{id})(R) = R^{1,3}R^{2,3}$,
implies that $R = R_1 R_2$.

Let us set\footnote{For $X$ a Hopf algebra, $X^{cop}$ means $X$ with
opposite coproduct; if $Y\subset X$, then $Y^{\varepsilon}
:= Y \cap \on{Ker}(\varepsilon_X)$.}
$B := A^{*cop}$,
$B_1 := (A_1^\varepsilon A_2)^\perp \subset B$,
$B_2 := (A_1 A_2^\varepsilon)^\perp\subset B$.

\begin{theorem} \label{thm:A1A2} (see \cite{ER,EF,DKhP})
\begin{itemize}
\item[1)] $B_1,B_2$ are subalgebras of $B$; the subalgebra
$B_1$ (resp., $B_2$) is a
left (resp., right) coideal of $B$, i.e., $\Delta_B(B_1) \subset
B \otimes B_1$, $\Delta_B(B_2) \subset B_2 \otimes B$, and $m_B: B_2 \otimes
B_1\to B$ is a vector space isomorphism.

\item[2)] Define $P'_i : B\to B_i$ by $P'_2(b_2b_1) = b_2 \varepsilon_B(b_1)$,
$P'_1(b_2b_1) = \varepsilon_B(b_2)b_1$. Then
$R_i = (\on{id}_A \otimes P'_{3-i})(R)$, for $i=1,2$.
In fact, $R_i = (P_i \otimes P'_{3-i})(R)\in A_i \otimes B_{3-i}$.

\item[3)] $R_2$ is a cocycle for $D$, i.e., $R_2^{1,2}
(\Delta_D \otimes \on{id}_D)(R_2) = R_2^{2,3}
(\on{id}_D \otimes \Delta_D)(R_2)$.
It follows that $D$, equipped with the coproduct ${}^{R_2}\Delta_D(x):=
R_2\Delta_D(x)R_2^{-1}$, is a quasitriangular Hopf algebra
(which we denote by ${}^{R_2}D$) with $R$-matrix $R_2^{2,1}R_1$.

\item[4)] $m_D(A_i \otimes B_i) = m_D(B_i \otimes A_i)$ for $i=1,2$, so
$D_i := m_D(A_i \otimes B_i) \subset D$ are subalgebras of $D$.

\item[5)] $A_i$, $B_i$ have the following coideal properties:
${}^{R_2}\Delta_D(A_1) \subset A_1 \otimes D_1$,
${}^{R_2}\Delta_D(B_1)
\subset D_1 \otimes B_1$, ${}^{R_2}\Delta_D(A_2) \subset A_2 \otimes D_2$,
${}^{R_2}\Delta_D(B_2) \subset D_2 \otimes B_2$.

\item[6)] $D_i$ are Hopf subalgebras of ${}^{R_2}D$. The quasitriangular
Hopf algebra ${}^{R_2}D$ is isomorphic to the double of
$(D_1,({}^{R_2}\Delta_D)_{|D_1})$, whose dual algebra with opposite
coproduct is $(D_2,({}^{R_2}\Delta_D)_{|D_2})$.
\end{itemize}
\end{theorem}

{\em Proof.} 1) For $X,Y\subset A$, set $XY: = m_A(X\otimes Y)$.

We have $(\on{id} \otimes \varepsilon_A) \circ \Delta_A = \on{id}_A$,
which implies that $\Delta_A(A_1^\varepsilon) \subset A \otimes
A_1^\varepsilon + A_1^\varepsilon \otimes 1_A$.

Then $\Delta_A(A_1^\varepsilon A_2) \subset \Delta_A(A_1^\varepsilon)
\Delta_A(A_2) \subset (A\otimes A_1^\varepsilon + A_1^\varepsilon \otimes 1_A)
(A_2\otimes A) \subset
A\otimes A_1^\varepsilon A + A_1^\varepsilon A_2\otimes A$.

Now $A = A_1A_2$, and $A_1^\varepsilon A_1= A_1^\varepsilon$, so
$A_1^\varepsilon A= A_1^\varepsilon A_2$, which
implies that $A_1^\varepsilon A_2$ is a two-sided coideal of $A$.
This implies that $B_1$ is a subalgebra of $B$. In the same way, $B_2$
is a subalgebra of $B$.

$A_1^\varepsilon A_2$ is a right ideal of $A$, which implies that $B_1$ is
a left coideal of $B$. In the same way, $B_2$ is a right coideal of $B$.

We now show that $m_B : B_2 \otimes B_1 \to B$ is a
vector space isomorphism. We have vector space isomorphisms
$B_i \simeq A_{3-i}^*$, for $i=1,2$, induced by $A = A_1 \oplus
A_1 A_2^\varepsilon$ and $A = A_2 \oplus A_1^\varepsilon A_2$.
So we will prove that the transposed map
$A \stackrel{m_B^t}{\to} A_1 \otimes A_2$ is an isomorphism.

Let $b_i\in A_i^* \simeq B_{3-i}$. When viewed as elements of $B = A^*$,
$b_i$ satisfy\footnote{We denote by $a\otimes b\mapsto \langle a,b\rangle =
\langle b,a\rangle$ the pairing $A\otimes B\to \CC$.}
$\langle b_1 , a_1a_2 \rangle = \langle b_1 , a_1 \rangle
\varepsilon(a_2)$, $\langle b_2 , a_1a_2 \rangle = \varepsilon(a_1)
\langle b_2 , a_2\rangle$. Then
$\langle b_1b_2 , a_1a_2 \rangle = \langle b_1 \otimes b_2,
a_1^{(1)}a_2^{(1)} \otimes a_1^{(2)}a_2^{(2)}\rangle
= \langle b_1 , a_1^{(1)} \rangle \varepsilon(a_2^{(1)})
\varepsilon(a_1^{(2)}) \langle b_2, a_2^{(2)} \rangle  =
\langle b_1 , a_1 \rangle \langle b_2, a_2 \rangle$.
So the composed map $A_1 \otimes A_2 \stackrel{m_A}{\to} A
\stackrel{m_B^t}{\to} A_1 \otimes A_2$ is the identity,
which implies that $A \stackrel{m_B^t}{\to} A_1 \otimes A_2$,
and therefore $m_B : B_2 \otimes B_1 \to B$, is an isomorphism.

2) Let $R'_i \in A_i \otimes B_{3-i}$ be the canonical element
arising from the isomorphism $B_i \simeq A_i^*$. Let us show that
$R = R'_1 R'_2$, and $R_i = R'_i$.

We have $R'_1R'_2\in A\otimes B$. For $a\in A$, let us compute
$\langle R'_1R'_2,\on{id} \otimes a \rangle$. We assume that $a = a_1a_2$,
with $a_i\in A_i$. Then
$\langle R'_1R'_2,\on{id} \otimes a \rangle
= \langle R'_1R'_2,\on{id} \otimes a_1a_2 \rangle
= \langle (R'_1)^{1,2}(R'_2)^{1,3},\on{id} \otimes a_1^{(1)}a_2^{(1)}
\otimes a_1^{(2)}a_2^{(2)} \rangle
= \langle (R'_1)^{1,2}(R'_2)^{1,3},\on{id} \otimes
a_1^{(1)}\varepsilon_A(a_2^{(1)})
\otimes \varepsilon_A(a_1^{(2)})a_2^{(2)} \rangle
= a_1a_2 = a$. So $R'_1R'_2\in A\otimes B$ is the canonical
element, so it is equal to $R$.

Now $R_i = (P_i\otimes \on{id}_B)(R) = R'_i$, since we have
$(\varepsilon_A \otimes \on{id}_B)(R_i) = 1_A$. We also compute
$(\on{id} \otimes P'_{3-i})(R) = R_i$ and $(P_i \otimes P'_{3-i})(R)
= R_i$.

3) We first prove that for any $a\in A$, $b\in B$, we have
\begin{equation}\label{D-mu}
\< a^{(1)}, b^{(1)}\> b^{(2)} a^{(2)}= a^{(1)}b^{(1)}
\< a^{(2)}, b^{(2)}\>\ .
\end{equation}
Recall the multiplication formula in the quantum double \cite{D}:
\begin{equation} \label{product:form}
b\ a = \< S_D(a^{(1)}),b^{(1)} \>\  \<a^{(3)},b^{(3)}\>\ a^{(2)}\ b^{(2)}\ .
\end{equation}
Using this equality, we may write the right hand side of \r{D-mu}
as follows:
\begin{align*}
&\< a^{(1)}, b^{(1)}\> \< S_D(a^{(2)}),b^{(2)} \>
\<a^{(4)},b^{(4)}\> a^{(3)} b^{(3)} =
\< a^{(1)}\ot S_D(a^{(2)}), \Delta_B(b^{(1)})\>  \<a^{(4)},b^{(3)}\>
a^{(3)} b^{(2)} \\
&\qquad = \varepsilon(a^{(1)})\varepsilon (b^{(1)})
\<a^{(3)},b^{(3)}\> a^{(2)}\ b^{(2)}=
\<a^{(2)},b^{(2)}\> a^{(1)}\ b^{(1)}
\end{align*}
which coincides with the left hand side of equality \r{D-mu}.

We now prove the cocycle relation $R_2^{1,2}R_2^{12,3}=
R_2^{2,3}R_2^{1,23}$ (the coproduct is $\Delta_D$). Both
sides of this identity belong to $A_2 \otimes D \otimes B_1$.
Using the pairing, we will identify them with linear maps
$B_1 \otimes A_2 \to D$.
Let us compute the pairing of both sides of this equality with
$b_1\otimes\on{id}\otimes a_2$ for arbitrary $b_1\in B_1$ and $a_2\in A_2$.
For the left hand side we have
\begin{align*}
& \< R_2^{1,2}R_2^{12,3},b_1\otimes\on{id}\otimes a_2\>=
\< R_2(a_2^{(1)}\ot a_2^{(2)}),b_1\otimes\on{id}\>=
\sum_i \< r'_i a_2^{(1)}\ot r''_i a_2^{(2)},b_1\otimes\on{id}\>\\
& \qquad = \sum_i \< r'_i\ot a_2^{(1)},b^{(2)}_1\otimes b^{(1)}_1\>
r''_i a_2^{(2)}=
\< a_2^{(1)}, b^{(1)}_1\> b^{(2)}_1 a_2^{(2)}\ .
\end{align*}
Here $R_2 = \sum_i r'_i \otimes r''_i$. On the other hand, for the right
hand side we obtain
\begin{align*}
& \< R_2^{2,3}R_2^{1,23},b_1\otimes\on{id}\otimes a_2\>=
\< R_2(b_1^{(1)}\ot b_1^{(2)}),\on{id}\ot a_2\>=
\sum_i \< r'_i b_1^{(1)}\ot r''_i b_1^{(2)}, \on{id} \ot a_2\>\\
& \qquad = \sum_i r'_i\ot b_1^{(1)} \< r''_i \ot b_1^{(2)},a^{(1)}_2
\otimes a^{(1)}_2\>=
a_2^{(1)}b^{(1)}_1 \< a_2^{(2)}, b^{(2)}_1\> \ ,
\end{align*}
where we used the fact that $R_2$ is the pairing tensor
between the subalgebras $A_2$ and $B_1$.
The cocycle identity now follows from (\ref{D-mu}).

Here is another proof of 3).
Recall that $R$ is invertible, and $R^{-1} = (S_D\otimes \on{id}_D)(R)$.
It follows that $R_1$ and $R_2$ are invertible. Let us show that
$R_1^{-1}\in A_1 \otimes B_2$. We first show that $R_1^{-1}\in D \otimes B_2$.
For this, we let $a_1\in A_1$, $a_2^\varepsilon  \in A_2^\varepsilon$
and we compute:
\begin{align*}
& \langle R_1^{-1}, \on{id} \otimes a_1a_2^\varepsilon \rangle =
\langle R_2 R^{-1}, \on{id} \otimes a_1a_2^\varepsilon \rangle =
\langle R_2^{1,2} (R^{-1})^{1,3}, \on{id} \otimes (a_1a_2^\varepsilon)^{(1)}
\otimes (a_1a_2^\varepsilon)^{(2)} \rangle
\\ &
= \langle R_2, \on{id} \otimes (a_1a_2^\varepsilon)^{(1)} \rangle
\langle R^{-1}, \on{id} \otimes (a_1a_2^\varepsilon)^{(2)}
\rangle
\\ &
= P_2((a_1a_2^\varepsilon)^{(1)})S_A((a_1a_2^\varepsilon)^{(2)})
= P_2(a_1^{(1)}) (a_2^\varepsilon)^{(1)} S_A((a_2^\varepsilon)^{(2)})
S_A(a_1^{(2)}) = 0,
\end{align*}
which proves that $R_1^{-1}\in D \otimes B_2$. In the same way, one proves
that $R_1^{-1}\in A_1 \otimes D$, so $R_1^{-1}\in A_1 \otimes B_2$,
and then $R_2^{-1}\in A_2 \otimes B_1$.

Let us set $\Phi := R_2^{2,3}R_2^{1,23}(R_2^{1,2}R_2^{12,3})^{-1}$
(the coproduct is $\Delta_D$). Then $\Phi\in A_2 \otimes D \otimes B_1$.
Using the quasitriangular identities satisfied
by $R$, we get $\Phi = (R_1^{1,\bar{23}} R_1^{2,3})^{-1}
R_1^{\bar{12},3} R_1^{1,2}$, where the coproduct is now
$\Delta_D^{2,1}(x) = R \Delta_D(x) R^{-1}$.
The last identity implies that $\Phi\in A_1 \otimes D \otimes
B_2$. Since $A_1 \cap A_2 = {\mathbb C} 1_A$ and $B_1 \cap B_2 =
{\mathbb C} 1_B$, we get $\Phi \in 1_D \otimes D \otimes 1_D$.
The pentagon identity satisfied by $\Phi$ then implies that $\Phi =
1_D^{\otimes 3}$.

4) Recall that for $a\in A$, $b\in B$, we have
$$
ab = \langle b^{(1)} , a^{(1)} \rangle  \langle b^{(3)} ,
S_A(a^{(3)}) \rangle b^{(2)}a^{(2)}, \quad b a = \langle b^{(1)} ,
S_A(a^{(1)}) \rangle
\langle b^{(3)} , a^{(3)} \rangle a^{(2)}b^{(2)}.
$$

Let us set $XY := m_D(X\otimes Y)$, for $X,Y\subset D$. We will
show that $B_1A_1 \subset A_1 B_1$.
Let $a_1\in A_1$, $b_1\in B_1$. Then $b_1 a_1
= \langle b_1^{(1)} , S_A(a_1^{(1)}) \rangle \langle b_1^{(3)},
a_1^{(3)} \rangle a_1^{(2)}b_1^{(2)}$; we have $\otimes_{i=1}^3
a_1^{(i)} \in A^{\otimes 2} \otimes A_1$ and $\otimes_{i=1}^3
b_1^{(i)} \in B^{\otimes 2} \otimes B_1$,
and since $\langle b,a \rangle = \varepsilon_B(b) \varepsilon_A(a)$
for $a\in A_1$, $b\in B_1$, we have
$b_1 a_1 = \langle b_1^{(1)} , S_A(a_1^{(1)}) \rangle a_1^{(2)}
b_1^{(2)}$;
now $b_1^{(1)} \otimes b_1^{(2)} \in B \otimes B_1$
and $a_1^{(1)} \otimes a_1^{(2)} \in A \otimes A_1$, which implies that
$b_1 a_1 \in A_1B_1$, as wanted. (One proves in the same way that $B_1 A_1
\subset A_1B_1$.)
It follows that $D_1 := A_1B_1$ is a subalgebra of $D$. In the same way,
one shows that $D_2 := A_2B_2$ is a subalgebra of $D$.

Note that as above, for $a_1\in A_1$, $b_1\in B_1$, we have
$$ a_1 b_1
= \langle b_1^{(1)} , a_1^{(1)} \rangle
\langle  b_1^{(3)} , S_A(a_1^{(3)}) \rangle b_1^{(2)}
a_1^{(2)} =  \langle b_1^{(1)} , a_1^{(1)} \rangle b_1^{(2)}
a_1^{(2)},
$$
 so that $B_1A_1 = A_1B_1\simeq A_1 \otimes B_1$.
In the same way, $B_2A_2 = A_2B_2 \simeq A_2 \otimes B_2$.

5) We have
${}^{R_2}\Delta_D(A_1) = R_2\Delta_D(A_1) R_2^{-1} \subset
R_2 (A\otimes A_1) R_2^{-1} \subset A\otimes D_1$, since
$R_2^{\pm 1}\in A_2 \otimes B_1$. On the other hand,
${}^{R_2}\Delta_D(A_1) = R_1^{-1}\Delta_D^{2,1}(A_1) R_1
\subset R_1^{-1} (A_1 \otimes A)R_1 \subset A_1 \otimes D$, since
$R_1^{\pm 1}\in A_1\otimes B_2$. Finally, ${}^{R_2}\Delta_D(A_1)
\subset A_1 \otimes D_1$. The other inclusions are proved similarly.

6). 4) and 5) imply that $D_i$ are Hopf subalgebras of $D$.

We have now $R_2^{2,1}R_1\in D_1 \otimes D_2$. It is
a nondegenerate tensor, as it is inverse to the pairing
$D_1 \otimes D_2 \simeq B_1A_1 \otimes A_2 B_2 \simeq
(A_1 \otimes B_1) \otimes (A_2 \otimes B_2)\to {\mathbb C}$,
given by the tensor product of the natural pairings
$A_1 \otimes B_2 \to {\mathbb C}$, $A_2 \otimes B_1
\to {\mathbb C}$.

Let us prove that $m_D : D_1 \otimes D_2 \to D$ is a vector space isomorphism.
The map $A_1 \otimes B_1 \otimes A_2 \otimes B_2 \simeq A_1 B_1 \otimes A_2B_2
= D_1 \otimes D_2 \to D \simeq A \otimes B$ is given by
$a_1 \otimes b_1 \otimes a_2 \otimes b_2 \mapsto
\langle b_1^{(1)} , a_2^{(1)} \rangle \langle b_1^{(3)} , a_2^{(3)}
\rangle a_1 a_2^{(2)} \otimes
b_1^{(2)} b_2$. One checks  that the inverse map is given by
$A_1 \otimes A_2 \otimes B_1 \otimes B_2 \simeq A \otimes B
\to D_1 \otimes D_2$ using the same formula, replacing $\Delta_D$
by ${}^{R_2}\Delta_D$.

The statement is now a consequence of the following fact: let
$(H,R_H)$ be a quasitriangular Hopf algebra, and $H_i$,
$i = 1,2$ be Hopf subalgebras, such that $R_H\in H_1\otimes H_2$ is
nondegenerate and $m_H : H_1 \otimes H_2 \to H$ is a vector space isomorphism,
then $H_2 = H_1^{*cop}$ and $H$ is the double of $H_1$ (indeed, since $R_H$
is nondegenerate, it sets up a vector space isomorphism $H_1 \simeq H_2^*$,
and since it satisfies the quasitriangularity equations, this is an
isomorphism $H_1 \simeq H_2^{*,cop}$ of Hopf algebras; we are then in the
situation of the theorem of \cite{D} on doubles).
\hfill \qed \medskip

\subsection{The graded case}

In the case when $A$ is a Hopf algebra in the category of
$\NN$-graded vector spaces with finite dimensional components,
the results of the previous section can be generalized as follows.

Let $(\alpha_{ij})_{1\leq i,j\leq r}$ be a nondegenerate matrix,
let $A'$ be a $\NN$-graded braided Hopf algebra, with finite dimensional
components and $A'[0] \simeq \CC$, where the braiding
is defined by $(q^{\alpha_{ij}})_{1\leq i,j\leq r}$. Let
$A := A' \otimes \CC[\ZZ^r]$ be the corresponding Hopf algebra.
Let $B'$ be the graded dual to $A'$ and $B$ be the corresponding
Hopf algebra. We then have a nondegenerate Hopf pairing $A\otimes B
\to \CC$. Let $D$ be the quotient of the bicrossproduct of $A$ and $B$
by the diagonal inclusion of $\CC[\ZZ^r]$.

%
%
%
%

To explain in what space $R$ lies, we introduce the following notion.
If $V = \oplus_{n\in{\mathbb Z}} V[n]$ is a ${\mathbb Z}$-graded
vector space, set $V^{\bar\otimes k} := \prod_{n_1,...,n_k\in{\mathbb Z}}
V[n_1]\otimes ... \otimes V[n_k]$, let $V^{\otimes_> k} \subset
V^{\bar\otimes k}$ be the
set of all combinations $\sum v_1 \otimes ... \otimes v_k$, such
that there exists a constant $c_1$ and functions $c_2(n_1)$, ...,
$c_k(n_1,...,n_{k-1})$, such that $\on{deg}(v_1) \geq c_1$,
$\on{deg}(v_{2}) \geq c_2(\on{deg}(v_1))$, ..., $\on{deg}(v_k)
\geq c_k(\on{deg}(v_1),...,\on{deg}(v_{k-1}))$.

Define $R_1\in A_1 \otimes_> B_2$ as the tensor of the pairing
$\langle-,-\rangle_1$, $R'_2\in A'_2 \otimes_> B'_1$ as the tensor
of $\langle-,-\rangle'_2$ and $R_0$ as the tensor of
$\langle -,- \rangle_0$. Then $R_0 = q^{r_0}$, where $r_0$ is inverse to the
matrix $(\alpha_{ij})$, and the tensor of the $\langle-,-\rangle_2$
is $R_2 := R'_2R_0$ (it belongs to a suitable extension of
$A_2 \otimes_> B_1$). The $R$-matrix of $D$ is then
$R = R_1R_2$ (also in a suitable extension of $A \otimes_> B$.


In Theorem 1, 1) is unchanged; 2) is unchanged, with the addition that
$B_i$ are now {\it graded} subalgebras of $B$; in 3), the cocycle identity
holds in $D^{\otimes_> 3}$ and the next statement is that
${}^{R_2}\Delta_D$ defines a topological bialgebra structure,
i.e., we have an algebra morphism
${}^{R_2}\Delta_D : D \to D^{\otimes_> 2}$ and coassociativity is an
identity of maps $D\to D^{\otimes_> 3}$;  4) is unchanged;
5) has to be understood in the topological sense; and 6) has to be
replaced by the statement that $D_i$ are topological subbialgebras of
${}^{R_2}D$.

\noindent
{\bf Example.}
One may take $A = U_q({\mathfrak b}_+)$, where ${\mathfrak b}_+$ is the
Borel subalgebra of a Kac-Moody Lie algebra, equipped with
the principal grading. In some cases, the twisted bialgebra
${}^{R_2}D$ is an ordinary Hopf algebra, i.e., ${}^{R_2}\Delta :
D\to D^{\otimes 2}$.

\subsection{The topological case} \label{assum}

If $V = \oplus_{n\in{\mathbb Z}} V[n]$ is a ${\mathbb Z}$-graded
vector space, let $V^{\otimes_< k}$ be the image of
$V^{\otimes_> k}$ by $v_1\otimes .. \otimes v_k
\mapsto v_k \otimes ... \otimes v_1$.
We define $(V^{\otimes_< k})[n]$ as the part of $V^{\otimes_< k}$
of total degree $n$ and $(V^{\otimes_< k})_{fs} :=
\oplus_{n\in\ZZ} (V^{\otimes_< k})[n]$ the `finite support'
part of $V^{\otimes_< k}$. We define $(V^{\otimes_> k})_{fs}$ similarly.
Then if $V$ is a ${\mathbb Z}$-graded algebra, we have algebra inclusions
$V^{\otimes k} \subset (V^{\otimes_< k})_{fs} \subset V^{\otimes_< k}$
and $V^{\otimes k} \subset (V^{\otimes_> k})_{fs} \subset V^{\otimes_> k}$.

We will make the following assumptions.

\begin{itemize}

\item[(H1)]
$D$ is a ${\mathbb Z}$-graded topological bialgebra.
Here topological means that the coproduct is an algebra morphism
$\Delta_D : D \to (D^{\otimes_< 2})_{fs}$ of degree $0$
(then the coassociativity is an equality of
maps $D \to (D^{\otimes_< 3})_{fs}$).

\item[(H2)]
$D\supset A,B \supset
\CC[\ZZ^r]$, where $A,B$ are $\ZZ$-graded topological subbialgebras
of $D$ and $\CC[\ZZ^r] \subset D_0$ is equipped with its standard
bialgebra structure. We assume that the product map
yields an isomorphism of vector spaces $A\otimes_{\CC[\ZZ^r]} B \to D$,
and that we have a nondegenerate bialgebra
pairing $\langle-,-\rangle : A\otimes B^{cop} \to {\mathbb C}$ of degree $0$;
this means that
\begin{equation} \label{bialg:pairing}
\langle a_1a_2,b\rangle = \langle a_1\otimes a_2, b^{(2)} \otimes b^{(1)}
\rangle, \quad
\langle a, b_1 b_2 \rangle = \langle a^{(1)}\otimes a^{(2)}, b_1 \otimes b_2
\rangle,
\end{equation}
$a_i\in A_i$, $b_i\in B_i$, and $\langle a, 1_D\rangle = \varepsilon_D(a)$,
$\langle 1_D, b\rangle = \varepsilon_D(b)$; we further require that
the identity (\ref{D-mu})
$$
\< a^{(1)}, b^{(1)}\> b^{(2)} a^{(2)}= a^{(1)}b^{(1)}
\< a^{(2)}, b^{(2)}\>
$$
holds. (One checks that all the sums involved in these identities are finite.)

\item[(H3)]
$A\supset A_1,A_2$, where $A_i$ are $\ZZ$-graded
subalgebras of $A$, such that
the nontrivial components of
$A_1$ (resp., $A_2$) are in degrees $\geq 0$ (resp., $\leq 0$),
and the product map
$A_1 \otimes A_2 \to A$ is a linear isomorphism.

$B\supset B_1,B_2$, where $B_i$ are $\ZZ$-graded subalgebras of $B$,
such that the nontrivial components of $B_1$ (resp., $B_2$) are
in degrees $\geq 0$ (resp., $\leq 0$) and the product map is an
isomorphism $B_2 \otimes B_1 \to B$.

\item[(H4)]
Moreover,
\begin{equation} \label{coideal:pties}
\Delta_D(A_1) \subset A\otimes_< A_1, \;
\Delta_D(A_2) \subset A_2 \otimes_< A, \;
\Delta_D(B_1) \subset B \otimes_< B_1, \;
\Delta_D(B_2) \subset B_2\otimes_< B.
\end{equation}
If we define $\langle-,-\rangle_i$ as the restriction of
$\langle-, -\rangle$ to $A_i\otimes B_{3-i}\to \CC$ ($i=1,2$), then
this assumption is equivalent to the identity
\begin{equation} \label{axiom:pairing}
\langle a_1a_2,b_2b_1\rangle
= \langle a_1,b_2\rangle_1 \langle a_2,b_1\rangle_2,
\end{equation}
where $a_i\in A_i$, $b_i\in B_i$ ($i=1,2$).

\item[(H5)]  The degree zero components are
$A_2[0] = B_1[0] = \CC[\ZZ^r]$,
and $A_1[0] = B_2[0] = \CC$. We also assume that $A_2$, $B_1$
contain graded subalgebras $A'_2$, $B'_1$, such that the product induces
linear isomorphims $A'_2 \otimes \CC[\ZZ^r] \simeq A_2$, $B'_1 \otimes
\CC[\ZZ^r] \simeq B_1$ (so $A'_2[0] = B'_1[0] = \CC$).
We assume that the homogeneous components of $A_1$, $A'_2$, $B'_1$, $B_2$
are finite dimensional.

\item[(H6)]
Let us denote by $\langle-, -\rangle_0$ the restriction of
$\langle-,-\rangle$ to $A_2[0] \otimes B_1[0] \simeq \CC[\ZZ^r]^{\otimes 2}
\to \CC$; we assume that it has the form
$\langle \delta_i, \delta_j \rangle_0 = q^{\alpha_{ij}}$, $q\in\CC^\times$
and the matrix $(\alpha_{ij})$ is nondegenerate (here $\delta_i$ is the $i$th
basis vector of $\ZZ^r$). Let us denote by $\langle-,-\rangle'_2$
the restriction of $\langle-,-\rangle_2$ to $A'_2\otimes B'_1\to\CC$.
We assume that
$$
\langle a'_2a_0, b'_1b_0\rangle_2 = \langle a'_2, b'_1\rangle'_2
\langle a_0,b_0\rangle_0
$$
where $a'_2\in A'_2$, $b'_1\in B'_2$, $a_0\in A_2[0]$ and $b_0\in B_1[0]$.
We assume that the pairings $\langle-,-\rangle_1$
and $\langle-,-\rangle'_2$ are non-degenerate, in the sense that
each pairing between each pair of finite-dimensional homogeneous
components of opposite degrees is nondegenerate.

\end{itemize}

Define $R_1\in A_1 \otimes_> B_2$ as the tensor of the pairing
$\langle-,-\rangle_1$, $R'_2\in A'_2 \otimes_< B'_1$ as the tensor
of $\langle-,-\rangle'_2$ and $R_0$ as the tensor of
$\langle -,- \rangle_0$. Then $R_0 = q^{r_0}$, where $r_0$ is inverse to the
matrix $(\alpha_{ij})$, and the tensor of the $\langle-,-\rangle_2$
is $R_2 := R'_2R_0\in A_2 \otimes_< B_1$.

Actually, $R_1,R_2$ have degree $0$, so we have
$R_1\in (A_1 \otimes_> B_2)_{fs}$,  $R_2\in (A_2 \otimes_< B_1)_{fs}$.
Since $R_1$ has the form $1 + \sum_{i>0}
a_i \otimes b_i$, where $\on{deg}(a_i) = -\on{deg}(b_i) = i$,
$R_1$ is invertible in $A_1 \otimes_> B_2$. In the same way, $R_2$ is
invertible.

\begin{lemma}
$R_2$ is a cocycle for $\Delta_D$, i.e., the identity
$$
R_2^{1,2}(\Delta_D \otimes \on{id})(R_2)
= R_2^{2,3}(\on{id} \otimes \Delta_D)(R_2)
$$
holds in $D^{\otimes_< 3}$.
In the same way, $R_1^{-1}$ is a cocycle for $\Delta_D^{2,1}$.
\end{lemma}

{\em Proof.} The proof for $R_2$ is the same as the first proof of 3)
in Theorem \ref{thm:A1A2}. In the case of $R_1^{-1}$, we similarly prove that
$(\Delta_D^{2,1} \otimes \on{id})(R_1) R_1^{1,2} = (\on{id} \otimes
\Delta_D^{2,1})(R_1) R_1^{2,3}$ and take inverses.
\hfill \qed \medskip

We therefore obtain topological bialgebra structures
$\Delta : D \to D^{\otimes_< 2}$ and $\tilde\Delta : D \to
D^{\otimes_> 2}$, defined by
$$
\Delta(x) := R_2 \Delta_D(x) R_2^{-1}, \quad
\tilde\Delta(x) = R_1^{-1} \Delta_D(x)^{2,1} R_1.
$$

We now prove:
\begin{theorem} \label{thm:twists}
$\Delta$ and $\tilde\Delta$ actually take their values in $D^{\otimes 2}$,
and are equal as maps $D\to D^{\otimes 2}$.
\end{theorem}

Let us first briefly summarize the proof. For $x\in A$, $\Delta(x)\in
A \otimes_< D$ while $\tilde\Delta(x)\in A \otimes_> D$. We pair both
elements with $b\otimes \on{id}$, where $b\in B$.
Then $\Delta(x)$, $\tilde\Delta(x)$ define elements of two completions
$\on{Hom}_{\pm}(B,\hat D)$ of the same convolution algebra
$\oplus_{(i,j)\in\ZZ^2} \on{Hom}(B_i,D_j)$. Using identities in these
convolution algebras, and computing degrees carefully to prove that
certain maps, a priori valued in $D^{\otimes_> 3}$ or
$D^{\otimes_< 3}$, take in fact their values in $D^{\otimes 3}$, we
prove identity (\ref{final}), which implies that the pairings of
$\Delta(x)$ and $\tilde\Delta(x)$ with $b\otimes \on{id}$ are the same.
This implies $\Delta(x) =\tilde\Delta(x)\in D^{\otimes 2}$; the proof
with $x$ replaced by $y\in B$ is similar.

\medskip

{\em Proof.} We will consider the
convolution algebra $\oplus_{(i,j)\in \ZZ^2}\on{Hom}(B_i,D_j)$,
where the product is
$(f_1 * f_2)(b) := f_1(b^{(2)}) f_2(b^{(1)})$. This is an
associative algebra with identity element
$1_* : b\mapsto \varepsilon(b)1_D$. This algebra is bigraded by
$\ZZ^2$. The convolution product can be extended as follows. Let
$\hat D := \prod_{i\in \ZZ} D_i$, then $\on{Hom}(B,\hat D) =
\prod_{(i,j)\in\ZZ^2} \on{Hom}(B_i,D_j)$.
If $f = \prod_{i,j} f_{i,j} \in \on{Hom}(B,\hat D)$,
we define the support of $f$ as $\on{supp}(f) :=
\{(i,j)|f_{i,j}\neq 0\}$. Then if
$f_1,f_2\in \on{Hom}(B,\hat D)$ are such that
the sum map $\on{supp}(f_1) \times \on{supp}(f_2) \to \on{supp}(f_1) +
\on{supp}(f_2)$ has finite fibers, then the convolution $f_1*f_2\in
\on{Hom}(B,\hat D)$ is defined, and has
support contained in $\on{supp}(f_1) + \on{supp}(f_2)$. One checks
that the convolution in $\on{Hom}(B,\hat D)$
is associative in the restricted sense that if
$S_1 \times S_2 \times S_3 \to S_1 + S_2 + S_3$ has finite fibers,
where $S_i := \on{supp}(f_i)$, then $(f_1*f_2)*f_3 = f_1 * (f_2 * f_3)$.

In particular, we define $\on{Hom}_+(B,\hat D)$ (respectively,
$\on{Hom}_-(B,\hat D)$) as the subset of
$\on{Hom}(B,\hat D)$ of all the elements
$f$ such that $\on{supp}(f)$ is contained is some part of $\ZZ^2$
of the form $S + \NN(1,1)$ (resp., $S + \NN(-1,-1)$),\footnote{$\NN(a,b) =
\{(0,0),(a,b),(2a,2b),...\}$}
where $S$ is a finite subset of $\ZZ^2$. Then $\on{Hom}_\pm(B,\hat D)$
are both algebras for the convolution.

One checks that one has algebra injections
$(A\otimes_< D)_{fs} \to \on{Hom}_+(B,\hat D)$ and $(A\otimes_> D)_{fs}
\to \on{Hom}_-(B,\hat D)$ (denoted $x\mapsto [x]$),
extending the map $A\otimes D\to
\on{Hom}(B,\hat D)$, $a\otimes d \mapsto
(b\mapsto \langle b,a\rangle d)$.
The intersection of $(A\otimes_< D)_{fs}$ and $(A\otimes_>D)_{fs}$ in
$\on{Hom}(B,\hat D)$ is $A\otimes D$.

It follows from (\ref{axiom:pairing}) that the images of
$R_1\in (A \otimes_> D)_{fs}$ and $R_2\in (A \otimes_< D)_{fs}$ in
$\on{Hom}(B,\hat D)$ respectively coincide with $P_2$ and $P_1$,
where $P_2(b_2b_1) = b_2 \varepsilon(b_1)$ and
$P_1(b_2b_1) = \varepsilon(b_2)b_1$.

Let $f:= [R_1^{-1}]\in \on{Hom}_-(B,\hat D)$ and $g:= [R_2^{-1}]\in
\on{Hom}_+(B,\hat D)$. Then $f*P_2$ and $P_2* f$ are defined, and
$f*P_2 = P_2 * f = 1_*$. This means that for any $b\in B$,
$f(b^{(2)})P_2(b^{(1)}) = P_2(b^{(2)})f(b^{(1)}) = \varepsilon(b)$.
For $b = b_2b_1$, this means that $f(b_2^{(2)}b_1^{(2)})b_2^{(1)}
P_2(b_1^{(1)})
= \varepsilon(b)$ and $P_2(b_2^{(2)})f(b_2^{(1)}b_1) = \varepsilon(b)$.
The first equality says in particular (setting $b_1 = 1$) that
$f(b_2^{(2)})P_2(b_2^{(2)}) =
\varepsilon(b_2)$; plugging this in the second equality, we get
$f(b_2) \varepsilon(b_1) = f(b_2^{(2)}) \varepsilon(b_2^{(1)}b_1) =
f(b_2^{(3)})P_2(b_2^{(2)})f(b_2^{(1)}b_1) = \varepsilon(b_2^{(2)})
f(b_2^{(1)}b_1) = f(b_2b_1)$, so
$$
f(b_2b_1) = f(b_2)\varepsilon(b_1).
$$
In addition, we have
$$
f(b_2^{(2)})b_2^{(1)} = P_2(b_2^{(2)})f(b_2^{(1)}) = \varepsilon(b_2).
$$
These identities imply
\begin{equation} \label{ids:f}
f(b^{(2)})b^{(1)} = P_1(b).
\end{equation}
In the same way, one shows that
$$
g(b_2b_1) = \varepsilon(b_2)g(b_1), \quad
g(b_1^{(2)})P_1(b_1^{(1)}) = b_1^{(2)}g(b_1^{(1)}) = \varepsilon(b_1),
$$
where all the a priori infinite sums reduce to finite sums.

Then we compute
\begin{align*}
& \langle \Delta(x), b\otimes \on{id}\rangle = (P_1 * [\Delta(x)] * g)(b)
= P_1(b^{(3)})[\Delta_D(x)](b^{(2)})g(b^{(1)})
\\ &
= P_1(b_2^{(3)}b_1^{(3)})\langle x^{(1)},b_2^{(2)}b_1^{(2)}\rangle
x^{(2)} g(b_2^{(1)}b_1^{(1)})
= P_1(b_2^{(2)})b_1^{(3)}\langle x^{(1)},b_2^{(1)}b_1^{(2)}\rangle
x^{(2)} g(b_1^{(1)})
\\ &
= P_1(b_2^{(2)})b_1^{(3)}\langle x^{(1)},b_2^{(1)}\rangle
\langle x^{(2)},b_1^{(2)}\rangle
x^{(3)} g(b_1^{(1)})
\end{align*}
where the fourth equality uses $b_1^{(3)}\in B_1$, $b_2^{(1)}\in B_2$, and the
last equality uses the bialgebra pairing rules
(\ref{bialg:pairing}), which do not introduce infinite sums.

We finally get
$$
\langle \Delta(x), b_2b_1\otimes \on{id}\rangle =
P_1(b_2^{(2)})\langle x^{(1)}, b_2^{(1)}\rangle \cdot
b_1^{(3)} \langle x^{(2)}, b_1^{(2)}\rangle x^{(3)} g(b_1^{(1)}),
$$
and similarly
$$
\langle \tilde\Delta(x), b_2b_1\otimes \on{id}\rangle =
f(b_2^{(3)}) \langle x^{(2)},b_2^{(2)} \rangle x^{(1)}b_2^{(1)} \cdot
\langle x^{(3)},b_1^{(2)} \rangle P_2(b_1^{(1)}).
$$

We now prove that
\begin{equation} \label{wanted:FG}
P_1(b_2^{(2)})\langle x^{(1)}, b_2^{(1)}\rangle  x^{(2)}=
f(b_2^{(3)}) \langle x^{(2)},b_2^{(2)} \rangle x^{(1)}b_2^{(1)} .
\end{equation}
We first show that both sides are finite sums. Let $x$ and $b_2$ be of
fixed degree (we denote by $|x|$ the degree of $x\in D$). Then for
some constant $c$, we have $|b_2^{(1)}|,|x^{(1)}| \leq c$,
$|b_2^{(2)}| = |b_2| - |b_2^{(1)}|$, $|x^{(2)}| = |x| - |x^{(1)}|$,
so the l.h.s. reduces to the sum of contributions with
$|b_2^{(1)}|,|x^{(1)}|\in \{-c,...,c-1,c\}$ and is
a finite sum. Let us show that the r.h.s. is a finite sum.
For some constant $c'$ and function $c''(n)$, we have
$|b_2^{(1)}|$, $|x^{(1)}|
\leq c'$, $|b_2^{(2)}|\leq c''(|b_2^{(1)}|)$,
$|x^{(2)}| = |x| - |x^{(1)}|$ and $|b_2^{(3)}| = |b_2| - |b_2^{(1)}|
- |b_2^{(2)}|$. The nontrivial contributions
are for $|x^{(2)}| + |b_2^{(2)}| = 0$ and $|b_2^{(3)}| \leq 0$
(as $\on{supp}(f) \subset \NN(-1,-1)$). These conditions
impose $|b_2^{(1)}| + |b_2^{(2)}| \geq 0$, so $|b_2^{(2)}| \geq -c'$,
which leaves only finitely many possibilities for $|b_2^{(2)}|$; then
$|b_2^{(1)}| \geq -|b_2^{(2)}|$, which leaves only finitely many
possibilities for $(|b_2^{(1)}|,|b_2^{(2)}|,
|b_2^{(3)}|)$. This also implies that only finitely many
$(|x^{(1)}|,|x^{(2)}|)$ contribute.

We have therefore proven that there are linear maps
$F,G : B_2 \otimes A \to D$,
$F(x\otimes b_2) := P_1(b_2^{(2)}) \langle x^{(1)},
b_2^{(1)} \rangle x^{(2)}$ and $G(x\otimes b_2) :=
f(b_2^{(3)}) \langle x^{(2)},b_2^{(2)} \rangle x^{(1)}b_2^{(1)}$.

Similarly, one proves that there exists a unique linear map
${\bold f} : B_2 \otimes A \to D$, such that
${\bold f}(b_2 \otimes x) = P_1(b_2^{(2)}) \langle x,b_2^{(1)} \rangle$.
Then the image of the composed map $B_2 \otimes A \stackrel{\on{id} \otimes
\Delta_D}{\to} B_2 \otimes A^{\otimes_< 2}
\stackrel{{\bold f} \otimes \on{inc}}{\to} D^{\otimes_< 2}$
(inc is the canoncial inclusion) is contained in $D^{\otimes 2}$, so its
composition with $m_D : D^{\otimes 2} \to D$ is well-defined.
Then $F = m_D \circ ({\bold f} \otimes \on{inc}) \circ (\on{id} \otimes
\Delta_D)$, and (\ref{wanted:FG}) expresses as $F = G$.

One checks that there are maps $u,v : A\otimes B \to D$, such that
$u(a\otimes b) := \langle a^{(1)},b^{(1)} \rangle b^{(2)} a^{(2)}$
and $v(a\otimes b) := \langle a^{(2)},b^{(2)} \rangle a^{(1)}b^{(1)}$;
(\ref{D-mu}) can then be expressed by the equality $u=v$.

As above, one checks that the composed map $A \otimes B \stackrel{\on{id}
\otimes \Delta_D}{\to} A \otimes B^{\otimes_< 2}
\stackrel{w\otimes \on{id}}{\to} D \otimes_< B \stackrel{\on{id}\otimes f}
{\to} D^{\otimes_< 2}$ actually takes its values in $D^{\otimes 2}$ for $w = u$
or $v$. This means that the equality
$$
\langle x^{(1)}, b_2^{(1)}\rangle b_2^{(2)}x^{(2)} \otimes f(b_2^{(3)})
= \langle x^{(2)}, b_2^{(2)} \rangle x^{(1)} b_2^{(1)} \otimes f(b_2^{(3)})
$$
takes place in $D^{\otimes 2}$. Applying $m_D$ after transposing the
factors, we get the identity in $D$
$$
f(b_2^{(3)})
\langle x^{(1)}, b_2^{(1)}\rangle b_2^{(2)}x^{(2)} = f(b_2^{(3)})
\langle x^{(2)}, b_2^{(2)} \rangle x^{(1)} b_2^{(1)},
$$
which according to (\ref{ids:f}) yields (\ref{wanted:FG}).

One proves similarly that the following is an equality between
finite sums
$$
b_1^{(3)} \langle x^{(1)}, b_1^{(2)}\rangle x^{(2)} g(b_1^{(1)})
= x^{(1)}\langle x^{(2)},b_1^{(2)} \rangle P_2(b_1^{(1)}),
$$
which one expresses as $F' = G'$, where $F',G' : A \otimes B \to D$
are given by $F'(x\otimes b_2):= b_1^{(3)} \langle x^{(1)}, b_1^{(2)}\rangle
x^{(2)} g(b_1^{(1)})$ and $G'(x\otimes b_2) :=
x^{(1)}\langle x^{(2)},b_1^{(2)} \rangle P_2(b_1^{(1)})$.

As before, there exists a unique linear map ${\bold g} : A \otimes B_1
\to D$, such that ${\bold g}(x\otimes b_1) = \langle x,b_1^{(2)}\rangle
P_2(b_1^{(1)})$. Then the image of the composed map $A\otimes B_1
\stackrel{\Delta_D \otimes \on{id}}{\to} A^{\otimes_< 2} \otimes B_1
\stackrel{\on{inc} \otimes {\bold g}}{\to} D^{\otimes_< 2}$ is contained in
$D^{\otimes 2}$, so its composition with $m_D : D^{\otimes 2} \to D$ is
well-defined, and $G' = m_D \circ (\on{inc} \otimes {\bold g}) \circ
(\Delta_D \otimes \on{id})$.

We then consider the composed map $B \otimes A\otimes B
\stackrel{\on{id} \otimes \Delta_D^{(2)} \otimes \on{id}}{\to}
B \otimes A^{\otimes_< 3} \otimes B \stackrel{{\bold f}\otimes \on{id}
\otimes {\bold g}}{\to} D^{\otimes_< 3}$ (where $\Delta_D^{(2)} =
(\Delta_D \otimes \on{id}) \circ \Delta_D$);
this map is
\begin{equation} \label{image}
b_2 \otimes x \otimes b_1 \mapsto P_1(b_2^{(2)}) \langle x^{(1)},b_2^{(1)}
\rangle  \otimes x^{(2)}\otimes \langle x^{(3)}, b_1^{(2)} \rangle
P_2(b_1^{(1)})
\end{equation}
and actually takes its values in $D^{\otimes 3}$, since for $|x|$ and $|b_2|$
fixed, we have $|b_1^{(1)}|,|b_2^{(1)}|,|x^{(1)}| \leq c$ for some $c$, and
the nontrivial contributions are for $|b_2^{(1)}| = -|x^{(1)}| \geq -c$,
which leaves only finitely many possibilities for $(|b_2^{(1)}|,|b_2^{(2)}|)$,
and therefore also for $|x^{(1)}|$. Now for each such $|x^{(1)}|$,
we have $|x^{(2)}| \leq c'(|x^{(1)}|)$ for some function $c'(n)$, and
$|x^{(3)}| = |x| - |x^{(1)}| - |x^{(2)}| \geq |x| - c - c'(|x^{(1)}|)$.
On the other hand, $|b_1^{(2)}| = |b_1| - |b_1^{(1)}| \geq |b_1| - c$,
Since we must have $|x^{(3)}| + |b_1^{(2)}| = 0$, this leaves only finitely
many possibilities for $(|b_1^{(2)}|,|x^{(3)}|)$. Finally, we have finitely
many possibilities for $(|b_1^{(1)}|,|b_1^{(2)}|)$ and for
$(|x^{(1)}|,|x^{(2)}|)$, hence for $(|x^{(1)}|,|x^{(2)}|,|x^{(3)}|)$.
So the r.h.s. of (\ref{image}) is a finite sum and belongs to $D^{\otimes 3}$.

We then consider the map
$$
m_D^{(2)} \circ ({\bold f} \otimes \on{id} \otimes {\bold g}) \circ
(\on{id} \otimes
\Delta_D^{(2)} \otimes \on{id})) : B \otimes A \otimes B \to D
$$
(where $m_D^{(2)} = (m_D \otimes \on{id}) \circ m_D$).
On one hand, we have
$$
m_D^{(2)} \circ ({\bold f} \otimes \on{id} \otimes {\bold g})
\circ (\on{id} \otimes  \Delta_D^{(2)} \otimes \on{id}))
= m_D \circ (F \otimes {\bold g}) \circ (\on{id} \otimes \Delta_D \otimes
\on{id})
= m_D \circ (G \otimes {\bold g}) \circ (\on{id} \otimes \Delta_D \otimes
\on{id});
$$
on the other hand, we have
$$
m_D^{(2)} \circ ({\bold f} \otimes \on{id} \otimes {\bold g})
\circ (\on{id} \otimes \Delta_D^{(2)} \otimes \on{id}))
= m_D \circ ({\bold f} \otimes G') \circ (\on{id} \otimes \Delta_D \otimes
\on{id})
= m_D \circ ({\bold f} \otimes F') \circ (\on{id} \otimes \Delta_D \otimes
\on{id});
$$
so
\begin{equation} \label{final}
m_D \circ (G \otimes {\bold g}) \circ (\on{id} \otimes \Delta_D \otimes
\on{id}) =
m_D \circ ({\bold f} \otimes F') \circ (\on{id} \otimes \Delta_D \otimes
\on{id}).
\end{equation}

Explicitly,
$$
P_1(b_2^{(2)}) \langle x^{(1)}, b_2^{(1)} \rangle \cdot \big(
x^{(2)}  \cdot \langle x^{(3)}, b_1^{(2)} \rangle P_2(b_1^{(1)}) \big) =
P_1(b_2^{(2)}) \langle x^{(1)}, b_2^{(1)} \rangle \cdot
b_1^{(3)} \langle x^{(2)}, b_1^{(2)} \rangle x^{(3)} g(b_1^{(1)})
$$
and
$$
\big( P_1(b_2^{(2)}) \langle x^{(1)}, b_2^{(1)} \rangle \cdot x^{(2)} \big)
\cdot \langle x^{(3)}, b_1^{(2)} \rangle P_2(b_1^{(1)}) =
f(b_2^{(3)}) \langle x^{(2)}, b_2^{(2)} \rangle x^{(1)} b_2^{(1)}
\cdot \langle x^{(3)}, b_1^{(2)} \rangle P_2(b_1^{(1)}),
$$
so (\ref{final}) is rewritten as
$$
P_1(b_2^{(2)})\langle x^{(1)}, b_2^{(1)}\rangle \cdot
b_1^{(3)} \langle x^{(2)}, b_1^{(2)}\rangle x^{(3)} g(b_1^{(1)})
= f(b_2^{(3)}) \langle x^{(2)},b_2^{(2)} \rangle x^{(1)}b_2^{(1)} \cdot
\langle x^{(3)},b_1^{(2)} \rangle P_2(b_1^{(1)}),
$$
i.e., $\langle \Delta(x), b_2b_1\otimes \on{id}\rangle =
\langle \tilde\Delta(x), b_2b_1\otimes \on{id}\rangle$.
This means that $[\Delta(x)] = [\tilde\Delta(x)]$. Since the
intersection of $(A\otimes_< D)_{fs}$ and $(A\otimes_> D)_{fs}$ in
$\on{Hom}(B,\hat D)$ is $A \otimes D$, we get
$\Delta(x) = \tilde\Delta(x) \in A\otimes D$.

One proves similarly that for $y\in B$, $\Delta(y) = \tilde\Delta(y)
\in D \otimes B$ by using the convolution algebra $\prod_{(i,j)\in\ZZ^2}
\on{Hom}(A_i,D_j)$,
where the product is given by $(f_1 * f_2)(a) = f_1(a^{(1)})f_2(a^{(2)})$.
\hfill \qed \medskip

\begin{proposition}
 $\Delta$ defines a (nontopological) bialgebra structure on $D$,
quasitriangular with $R$-matrix $R_2^{2,1}R_1\in D^{\otimes_> 2}$
(the quasitriangular identities are satisfied in $D^{\otimes_> 3}$).
\end{proposition}

{\em Proof.} Let us prove that for $x\in D$,
\begin{equation} \label{RDelta}
R_2^{2,1}R_1 \Delta(x) = \Delta^{2,1}(x) R_2^{2,1}R_1
\end{equation}
(equality in $D^{\otimes_>2}$).
We have $\Delta(x) = \tilde\Delta(x) = R_1^{-1} \Delta_D^{2,1}(x) R_1$
($\in D^{\otimes_>2}$) so the l.h.s. is equal to $R_2^{2,1}
\Delta_D^{2,1}(x)R_1$ ($\in D^{\otimes_> 2}$).

On the other hand, $\Delta^{2,1}(x) = R_2^{2,1} \Delta_D^{2,1}(x)
(R_2^{2,1})^{-1}$ ($\in D^{\otimes_> 2}$) so the r.h.s. is equal to
$R_2^{2,1} \Delta_D^{2,1}(x) R_1$ ($\in D^{\otimes_>2}$).
This proves (\ref{RDelta}).

We now prove the quasitriangular identity
\begin{equation} \label{qt:R}
(\Delta \otimes \on{id})(R) = R^{1,3} R^{2,3}
\end{equation}
in $D^{\otimes_> 3}$. Recall that
$R_2^{1,2}(\Delta_D\otimes\on{id})(R_2) = R_2^{2,3}(\on{id} \otimes
\Delta_D)(R_2)$ (equality in $D^{\otimes_< 3}$), which gives by applying
the transposition $x\mapsto x^{3,2,1}$
\begin{equation} \label{cha1}
R_2^{3,2}(\on{id} \otimes \Delta_D^{2,1})(R_2^{2,1}) =
R_2^{2,1}(\Delta_D^{2,1}\otimes \on{id})(R_2^{2,1})
\end{equation}
(equality in $D^{\otimes_> 3}$). On the other hand, recall that
\begin{equation} \label{cha2}
(\Delta_D^{2,1} \otimes \on{id})(R_1) R_1^{1,2} = (\on{id} \otimes
\Delta_D^{2,1})(R_1) R_1^{2,3}
\end{equation}
(equality in $D^{\otimes_> 3}$). Taking the product of (\ref{cha1})
written in opposite order
with (\ref{cha2}), we get
$$
R_2^{2,1}(\Delta_D^{2,1}\otimes \on{id})(R) R_1^{1,2}
= R_2^{3,2} (\on{id} \otimes \Delta_D^{2,1})(R) R_1^{2,3}.
$$
Equality (\ref{qt:R}) then follows from $R_1^{1,2}\Delta(x)
= \Delta_D^{2,1}(x) R_1^{1,2}$ (equality in $D^{\otimes_> 2}$).
The proof of the identity $(\on{id}\otimes\Delta)(R) = R^{1,3}R^{1,2}$
is similar.
\hfill \qed \medskip

\begin{proposition} \label{prop:coid:gal}
$D_1 := A_1 B_1$, $D_2 := A_2 B_2$ are subbialgebras of $D$.
We also have
\begin{equation} \label{inc:delta}
\Delta(A_1) \subset A_1 \otimes D_1, \quad
\Delta(B_1) \subset D_1 \otimes B_1, \quad
\Delta(A_2) \subset A_2 \otimes D_2, \quad
\Delta(B_2) \subset D_2 \otimes B_2.
\end{equation}
\end{proposition}

{\em Proof.} We first prove that $B_1A_1 \subset A_1B_1$.
Set $R(a,b) := \langle a^{(1)},b^{(1)} \rangle b^{(2)}a^{(2)} -
\langle a^{(2)},b^{(2)} \rangle a^{(1)}b^{(1)}$ and
$S(a,b) := ba - \langle S_D(a^{(1)}),b^{(2)}\rangle
\langle a^{(3)},b^{(3)}\rangle a^{(2)}b^{(2)}$.
We have $R(a,b) = \langle a^{(1)},b^{(1)}\rangle S(a^{(2)},b^{(2)})$
and $S(a,b) = \langle a^{(1)},S_D(b^{(1)}) \rangle R(a^{(2)},b^{(2)})$
so since $R(a,b) = 0$, we get $S(a,b) = 0$. Therefore
$$
ba = \langle S_D(a^{(1)}),b^{(2)}\rangle
\langle a^{(3)},b^{(3)}\rangle a^{(2)}b^{(2)}.
$$
If now $a\in A_1$, $b\in B_1$, we get $a^{(1)} \otimes a^{(2)} \otimes a^{(3)}
\in A^{\otimes_< 2} \otimes_< A_1$ and $b^{(1)} \otimes b^{(2)} \otimes b^{(3)}
\in B^{\otimes_< 2} \otimes_< B_1$. For $x\in A_1$, $y\in B_1$, we have
$\langle x,y \rangle = \varepsilon(x)\varepsilon(y)$, so
$ba = \langle S_D(a^{(1)}),b^{(1)}\rangle
\varepsilon(a^{(3)}) \varepsilon(b^{(3)})a^{(2)}b^{(2)}
= \langle S_D(a^{(1)}),b^{(1)}\rangle
a^{(2)}b^{(2)}$. Now since $\Delta_D(a_1)\in A \otimes_< A_1$ and
$\Delta_D(b_1)\in B \otimes_< B_1$, we get $ba\in A_1 B_1$, as wanted.
This implies that $D_1 := A_1B_1$ is a subalgebra of $D$.

In the same way, we prove that $D_2$ is a subalgebra of $D$.

 Let us prove that
$\Delta(A_1) \subset A_1 \otimes D_1$. For $x\in A_1$,
$\Delta(x) = R_2 \Delta_D(x) R_2^{-1}$ and since $D_1$ is an algebra,
$\Delta(x)\in A \otimes_< D_1$. On the other hand,
$\tilde\Delta(x) = R_1^{-1} \Delta_D^{2,1}(x)R_1 \in A_1 \otimes_> D$.
Since $\Delta(x) = \tilde\Delta(x)\in A \otimes D$ and $(A\otimes D)
\cap (A \otimes_< D_1) = A \otimes D_1$, we get $\Delta(x)\in
A \otimes D_1$. In the same way, $\Delta(x) \in A_1\otimes D$.
Then $(A\otimes D_1) \cap (A_1 \otimes D) = A_1 \otimes D_1$,
which implies $\Delta(x)\in A_1 \otimes D_1$, as wanted.
The other inclusions (\ref{inc:delta}) are proved in the same
way.

Since $D_i$ are generated by $A_i$, $B_i$, these inclusions
imply that $D_i$ are also subbialgebras of $(D,\Delta)$.
\hfill \qed \medskip

We will show that the quantum affine algebras, equipped with their
currents coproducts, are examples of the situation of
Subsection \ref{assum}.

\section{Quantum affine algebra $\Uqg$}\label{sl3-intro} \label{sect:doubles}

In this paper, $q$ is a complex number, which is neither $0$ nor a root
of unity.

\subsection{Chevalley-type presentation of $\Uqg$}\label{ch-gen}

Let $\ggg$ be a simple Lie algebra; let $r$ be its rank and let
$(b_{i,j})_{i,j=1,...,r}$ be its Cartan matrix. Let
$(a_{i,j})_{i,j=0,...,r}$ be the Cartan matrix of the affine Lie
algebra $\widehat{\ggg}$.
We denote by $\Pi=\{\a_1,...,\a_r\}$ the set of
positive simple roots of $\ggg$ and by $\widehat\Pi=\{\a_0,\a_1,...,\a_r\}$
the set of positive simple roots of $\widehat{\ggg}$.
The symmetrized Cartan matrix of $\widehat{\ggg}$ is
$((\a_i,\a_j))_{i,j = 0,...,r}$; we have
$(\a_i,\a_j)=d_i a_{i,j}= d_j a_{j,i}$ (where $d_i = 1,2$ or $3$
are coprime). Let $\delta$ be the minimal positive imaginary
root of $\widehat{\ggg}$, so
$\delta=\sum_{i=0}^r n_i\a_i$, $n_i\in \ZZ_{\geq 0}$, $n_0=1$. Let
$\left[{\gfr{n}{k}}\right]_q=\frac{[n]_q!}{[k]_q![n-k]_q!}$,
$[n]_q!=[1]_q[2]_q\cdots [n]_q$, $[n]_q=\frac{q^n-q^{-n}}{q-q^{-1}}$,
 $q_\a=q^{\frac{(\a,\a)}{2}}$, $q_i=q_{\a_i}=q^{d_i}$.

The quantum (untwisted) affine Lie algebra $\Uqg$ is generated by the
Chevalley generators $e_{\pm\a_i}$, $k^{\pm1}_{\a_i}$ ($i=0,\ldots,r$),
the grading elements $q^{\pm d}$, and the central elements
$k_\delta^{\pm 1/2}$, subject to the relations
$$
[q^d, k_{\alpha_i}] = [k_{\alpha_i},k_{\alpha_j}]=0,
\; q^d e_{\pm \a_i}q^{-d} =
q^{\pm\d_{i,0}}e_{\pm \a_i}, \;
k_{\al_i}e_{\pm\al_j}k^{-1}_{\al_i}\,=\,q_i^{\pm a_{ij}}e_{\pm\al_j},
$$
$$
(k_\delta^{\pm 1/2})^2 = \prod_{i=0}^r k_{\alpha_i}^{\pm n_i},  \quad
q^dq^{-d} = q^{-d}q^d = k_{\alpha_i}k_{\alpha_i}^{-1} = k_{\alpha_i}^{-1}
k_{\alpha_i} = k_\delta^{1/2} k_\delta^{-1/2} = k_\delta^{-1/2} k_\delta^{1/2}
= 1,
$$
\begin{equation}
[e_{\al_i},e_{-\al_j}]\,=\,
\delta_{ij}\frac{k_{\al_i}-k^{-1}_{\al_i}}{q_i-q_i^{-1}}\,,
\label{Q1}
\end{equation}
$$\sum_{r+s=1-a_{i,j}}(-1)^r e_{\pm\a_i}^{(r)}e_{\pm\a_j} e_{\pm\a_i}^{(s)}
=0,
\quad i\not= j\,,
\qquad
\text{where}\quad e_{\pm\a_i}^{(k)}=
\dfrac{e_{\pm\a_i}^{k}}{[k]_{q_i}!}\,.$$

The {\it standard Hopf structure} of $\Uqg$ is given by the formulas:
\begin{equation}
\begin{split}
\Delta^{std}(q^{\pm d})&=q^{\pm d}\ot q^{\pm d}, \qquad
\Delta^{std}(k_{\al_i}^{\pm 1})\,=\,k_{\al_i}^{\pm 1}\otimes k_{\al_i}^{\pm 1},
\qquad
\Delta^{std}(k_\delta^{\pm 1/2}) =
k_\delta^{\pm 1/2} \otimes k_\delta^{\pm 1/2},
\\
\Delta^{std}(e_{\al_i})&=e_{\al_i}\otimes 1+k_{\al_i}\otimes e_{\al_i}\,,\quad
\Delta^{std}(e_{-\al_i})=1\otimes e_{-\al_i}+e_{-\al_i}\otimes
k^{-1}_{\al_i},
\\
\coun(q^{\pm d})&=1,\qquad
\coun(e_{\pm\al_i})=0\,,\quad  \coun(k^{\pm1}_{\al_i})=1\,, \quad
\coun(k_\delta^{\pm 1/2}) = 1,
 \label{copr}\\
\ant(e_{\al_i})&=-k^{-1}_{\al_i}e_{\al_i}, \;
\ant(e_{-\al_i})=-e_{-\al_i}k_{\al_i}, \;
\ant(k^{\pm1}_{\al_i}) =  k^{\mp1}_{\al_i},
\;\\
\ant(q^{\pm d})&=q^{\mp d}, \;
\ant(k_\delta^{\pm 1/2}) = k_\delta^{\mp 1/2},
\end{split}
\end{equation}
where $\Delta^{std}$, $\coun$ and $\ant$ are the coproduct, counit and
antipode maps respectively.

Let $U_q(\h)$ be the Cartan subalgebra of  $\Uqtri$. It is  generated by the
elements $k_{\a_i}^{\pm 1}$ ($i = 0,...,r$) and $q^{\pm d}$.
Denote by  $U_q(\mathfrak{b}_+)$ the subalgebra of $\Uqg$ generated by
the elements $e_{\a_i}$, $k_{\a_i}^{\pm 1}$ ($i = 0,...,r$),
$k_\delta^{\pm 1/2}$ and $q^{\pm d}$,
and by $U_q(\mathfrak{b}_-)$ the subalgebra of $\Uqg$ generated by
the elements $e_{-\a_i}$, $k_{\a_i}^{\pm 1}$ ($i = 0,...,r$),
$k_\delta^{\pm 1/2}$ and $q^{\pm d}$.

The algebras $U_q(\mathfrak{b}_\pm)$ are Hopf subalgebras of $\Uqtri$
with respect to the standard coproduct $\Delta^{std}$. They are
$q$-deformations of the enveloping algebras of opposite Borel subalgebras of
Lie algebra $\widehat{\mathfrak{g}}$. We call them the {\it standard Borel
subalgebras}. Moreover, $U_q({\mathfrak b}_-)$ is the dual, with opposite
coproduct, of $U_q({\mathfrak b}_+)$, and $U_q(\widehat{\mathfrak g})
\otimes U_q({\mathfrak h})$ is the double of $U_q({\mathfrak b}_+)$
(where $U_q({\mathfrak h})$ is equipped with the standard structure,
for which $k_{\alpha_i}^{\pm 1}$ is primitive)

The algebras $U_q({\mathfrak b}_\pm)$ contain subalgebras
$U_q(\mathfrak{n}_\pm)$, which are generated by the elements
$e_{\pm \a_i}$, $i = 0,...,r$.
The subalgebra $U_q(\mathfrak{n}_+)$ is a left coideal of
$U_q(\mathfrak{b}_+)$ with respect to standard coproduct and the
subalgebra $U_q(\mathfrak{n}_-)$ is a right coideal of
$U_q(\mathfrak{b}_-)$ with respect to the same coproduct, that is
$$
\Delta^{std}(U_q(\mathfrak{n}_+))\subset U_q(\mathfrak{b}_+)\otimes
U_q(\mathfrak{n}_+)\,, \qquad
\Delta^{std}(U_q(\mathfrak{n}_-))\subset U_q(\mathfrak{n}_-)\otimes
U_q(\mathfrak{b}_-) \, .
$$
The algebras $U_q(\mathfrak{n}_\pm)$ are $q$-deformations of the
enveloping algebras of the standard pro-nilpotent subalgebras of
the affine Lie algebra $\widehat{\mathfrak{g}}$.

\subsection{Cartan-Weyl basis of $U_q(\widehat\ggg)$}
\label{sectionCW}

We now construct a Cartan-Weyl (CW) basis of $\Uqg$.
Let us denote by $\Sigma_+ \subset \Sigma$ the set of positive (resp.,
all) roots of $\ggg$ and by $\widehat\Sigma_+ \subset \widehat\Sigma$
the set of positive (resp., all) roots of $\widehat\ggg$, so
$\Sigma_+ \subset \Sigma$, $\widehat\Sigma_+ \subset \widehat\Sigma$.
Recall that $\a_0$ is the affine positive simple root and let
$\delta$ be the minimal positive imaginary root, so $\delta=\a_0+\theta$,
where $\theta$ is the longest root of $\Sigma_+$.

We consider the following class of normal orderings on $\widehat\Sigma_+$.
Let $\widehat{W}$ be the Weyl group of $\widehat{\ggg}$. It contains the
Weyl group $W$ of $\ggg$ and the normal subgroup $Q$, which is the set of all
elements having only finitely many conjugates. There is a unique group
morphism $Q\to \h^*$, $p\mapsto \bar p$, such that
the action of $p\in Q$ on $\h^*$ is the translation by $\bar{p}$.
This map is injective and identifies $Q$ with a subgroup
$\bar{Q}$ of $\h^*$.

Choose $p\in Q$ such that $(\bar{p},\a_i)>0$ for any $i=1,...,r$.
Choose a reduced decomposition
$p=s_{\a_{i_0}}s_{\a_{i_1}}\cdots s_{\a_{i_{m-1}}}$,
such that $\a_{i_0}=\a_{0}$ is the affine positive root. Extend the sequence
$i_0,i_1,...,i_m$ to a periodic sequence
\begin{equation}\label{period}
 \ldots, i_{-1},i_0,i_1,\ldots,i_n,\ldots
\end{equation}
satisfying the conditions $i_n=i_{n+m}$ for any $n\in\ZZ$.
We then set
$$
\gamma_1 := \alpha_{i_1}, \;
\gamma_2 := s_{\alpha_{i_1}}(\alpha_{i_2}), \; ..., \;
\gamma_k := s_{\alpha_{i_1}}...s_{\alpha_{i_{k-1}}}(\alpha_{i_k}) \;
\on{for} \;  k\geq 1;
$$
and
\begin{equation} \label{gammak}
\gamma_0 := \alpha_{i_0}, \;
\gamma_{-1} := s_{\alpha_{i_0}}(\alpha_{i_{-1}}), \;  ..., \;
\gamma_{-\ell} := s_{\alpha_{i_0}}...s_{\alpha_{i_{1-\ell}}}
(\alpha_{i_{-\ell}}) \; \on{for} \; \ell\geq 0.
\end{equation}


Then \cite{Be,Dam} the order $\gamma_1\prec\gamma_2\prec...\prec
\gamma_n\prec...\prec \delta
\prec 2\delta\prec...\prec\gamma_{-n}\prec...\prec
\gamma_{-1}\prec\gamma_0$ is normal and satisfies the condition
\begin{equation}\label{ord1}
l\delta+\a \prec (m+1)\delta \prec (n+1)\delta-\beta\,,
\end{equation}
for any positive roots $\a,\beta\in \Sigma_+$, and any $l,m,n\geq 0$.
From now on, we fix a normal ordering $\ord$ on $\widehat\Sigma_+$, given
by the procedure above.

Recall that the principal degree deg is the linear additive map
$\ZZ\widehat\Pi \to \ZZ$, such that $\operatorname{deg}(\alpha_i) = 1$
for $i = 0,...,r$. We first construct the CW generators $e_\gamma$,
for $\gamma\in \widehat\Sigma_+$ and $\operatorname{deg}(\gamma) \leq
\operatorname{deg}(\delta)-1$. For $\gamma\in \widehat\Pi$, $e_\gamma$
is equal to the corresponding Chevalley generator of $U_q(\widehat\ggg)$.
The CW generator $e_\gamma$ is then constructed by induction on the degree of
$\gamma$ as follows: if $\gamma\in
\widehat\Sigma_+$ and $\on{deg}(\gamma) \leq \on{deg}(\delta)-1$,
we let $[\alpha,\beta]$ be minimal for the inclusion, among the
set of all segments $[\alpha',\beta']$, where $\alpha',\beta'\in
\widehat\Sigma_+$ are such that $\gamma = \alpha' + \beta'$
(we define the segment $[\alpha',\beta']$ as $\{\gamma' |
\alpha' \preceq \gamma'\preceq \beta'\}$). We then set
\begin{equation}\begin{split}
\label{CW1}
e_\ga& :=[e_\a,e_\b]_{q^{-1}}\,,\qquad
e_{-\ga} :=
[e_{-\b},e_{-\a}]_{q^{}},
\end{split}\end{equation}
where\footnote{One can introduce the algebra
$\tilde U_q(\widehat{\mathfrak g})$ over the ring
$\CC[\tilde q,\tilde q^{-1},1/(\tilde q^n-1); n\geq 1]$
with the same generators and relations as $U_q(\widehat{\mathfrak g})$;
one can show using the CW basis that this is a free module over this
ring, and $U_q({\mathfrak g})$ is its specialization.
$\tilde U_q(\widehat{\mathfrak g})$ is equipped with the Cartan
antiinvolution $x\mapsto x^*$, defined by
$e_{\pm\a_i}^*=e_{\mp\a_i}$, $k_{\a_i}^*=k_{\a_i}^{-1}$, $\tilde q^*
=\tilde q^{-1}$. Then the analogues of $e_{\pm \gamma}$ satisfy
$e_{-\gamma}=e_\gamma^*$} the $q$-commutator $[e_\a,e_\b]_q$ means
$[e_\a,e_\b]_q=e_\a e_\b-q^{(\a,\b)}e_\b e_\a$.

The basis elements defined in this way coincide (up to normalization)
with those defined by the braid group action, using the same ordering
on roots. Since the latter basis is convex, $e_\gamma$ defined
above is independent (up to normalization) of the choice of the segment
$[\alpha,\beta]$, and depends only on the ordering of roots. Moreover,
it has been shown in \cite{Dam}, Proposition 11 that $e_{\delta - \alpha_i}$
is also independent on the choice of a normal ordering (up to normalization).

We then put
$$
e_\d^{(i)}=e_\d^{\prime(i)}=[e_{\a_i},e_{\d-\a_i}]_{q^{-1}},
\qquad e_{-\d}^{(i)}= [e_{\alpha_i - \delta},e_{-\alpha_i}]_q
$$
and\footnote{As before, the analogues of $e_{\pm\delta}^{(i)}$
in $\tilde U_q(\widehat{\mathfrak g})$ satisfy $e_{-\delta}^{(i)} =
(e_\delta^{(i)})^*$.} by induction for all $k>0$
\begin{equation*}\label{im2}
e_{\pm(\a_i+k\d)}=\pm\frac{1}{[2]_{q_i}}
[e_{\pm(\a_i+(k-1)\d)},e_{\pm\d}^{(i)}],\quad
e_{\pm(\d-\a_i+k\d)}=\pm\frac{1}{[2]_{q_i}}
[e_{\pm\d}^{(i)},e_{\pm(\d-\a_i+(k-1)\d)}]\ ,
\end{equation*}
\begin{equation*}\label{im4}
e^{\prime(i)}_{k\d}=
[e_{\a_i+(k-1)\d},e_{\d-\a_i}]_{q^{-1}}\ ,
\quad
e^{\prime(i)}_{-k\d}=
[e_{-\d+\a_i},e_{-\a_i-(k-1)\d}]_{q^{}}
\end{equation*}
\begin{equation}\label{im-rel}
e_{\pm k\d}^{(i)}=\!
\sum_{p_1+2p_2+\cdots+np_n=k}
\frac{(q_i^{\mp1}-q_i^{\pm 1})^{\sum p_j-1}
(\sum p_j-1)!}{p_1!\cdots p_n!}
\left(e^{\prime(i)}_{\pm\d}\right)^{p_1}\cdots
\left(e^{\prime(i)}_{\pm n\d}\right)^{p_n}\!\!\!\!.
\end{equation}
Then we apply procedure \rf{CW1} again to construct the remaining
(real root) CW
generators. As before, if $\gamma$ is real, then $e_\gamma$ depends
only on the choice of a normal ordering (up to normalization); the
$e_{n\delta \pm \alpha_i}$ and the $e^{(i)}_{n\delta}$ are independent
on the choice of this ordering (up to normalization).

The CW generators of the Borel subalgebras $U_q(\mathfrak{b}_\pm)$ satisfy
the following properties (see \cite{KhT2,Be,Dam})
\begin{align}
k_\a e_{\b} k_\a^{-1}&=q^{(\a,\b)}e_{\b}, \qquad
[e_\a,e_{-\a}]=a(\a)\frac{k_\a-k_\a^{-1}}{q-q^{-1}},
\notag \\
[e_{\pm\a},e_{\pm\b}]_{q^{-1}}&=\sum_{\{\ga_j\},\{n_j\}}
 C_{\{n_j\}}^{\pm\{\ga_j\}}(q)\
e^{n_1}_{\pm\ga_1}e^{n_2}_{\pm\ga_2}\cdots e^{n_m}_{\pm\ga_m},
\qquad \a,\b\in\Gama \label{use-pro}
\end{align}
where in \r{use-pro}  the sum is over all $\ga_1,\ldots,\ga_m\in \Sigma_+$,
$n_1,\ldots,n_m>0$ such that
$\a\ord\ga_1\ord\ga_2\ord\cdots\ord\ga_m\ord\b$ and
$\sum_j n_j\ga_j=\a+\b$ and the coefficients
$C_{\{n_j\}}^{\{\ga_j\}}(q)$ and $a(\a)$ are rational functions in
$q$ in ${\mathbb Q}[q,q^{-1},1/(q^n-1);n\geq 1]$.
The elements $k_\a$, $\a\in\widehat\Sigma$ are defined according to the
prescriptions $k_{\a+\b}=k_\a k_\b$, $k_\a=k_{\a_i}^{\pm 1}$, if
$\a=\pm\a_i$, $i = 0,...,r$.
The analogues of $e_{\pm\alpha}$ in $\tilde U_q(\widehat{\mathfrak g})$
also satisfy $e_{-\alpha} = e_\alpha^*$

The commutators $[e_\a,e_{-\b}]$, and $[e_{-\a},e_{\b}]$,
where $\a,\b\in\widehat\Sigma_+$, satisfy properties analogous to \rf{use-pro},
but the structure coefficients $C_{\{n_j\}}^{\pm\{\ga_j\}}$ belong to
$U_q({\mathfrak h})$. One can construct slightly
different generators in the CW basis such that property \r{use-pro}
is still valid with scalar structure constants. For this,
the normal ordering $\prec$ in the system $\widehat\Gama$
of positive roots is extended to a `circular' normal ordering (which is
not a total order) of the system $\Gam$ of all roots of $\widehat{\ggg}$.
This order $\prec_c$ is defined by: for $\alpha,\beta\in \Sigma_+$,
$(\alpha\prec_c \beta)\Leftrightarrow (\alpha\prec\beta) \Leftrightarrow
(-\alpha\prec_c -\beta)$ and
$(\alpha\prec_c -\beta)\Leftrightarrow (\beta\prec\alpha) \Leftrightarrow
(-\alpha\prec_c \beta)$. The order $\prec_c$ can be described as follows.
Suppose that $\gamma_1\prec\gamma_2\prec...\prec\gamma_N\prec...$
is the normal ordering of $\Gama$. Then we put all  the roots of $\Gam$ on
the circle clockwise in the following order (see Fig.~\ref{fig1}):
\begin{equation}\label{circle}
\gamma_1,\gamma_2,\ldots,\gamma_N,\ldots,
-\gamma_1,-\gamma_2,\ldots,-\gamma_N,\ldots,
\end{equation}
and say that the root $\gamma'\in \Gam$ precedes the root
$\gamma''\in \Gam$, $\gamma'\prec_c \gamma''$, if the segment
$[\gamma',\gamma'']$ in the circle \rf{circle} does not contain any of
the opposite roots $-\gamma'$ or $-\gamma''$.

For any positive root $\gamma\in\widehat\Sigma_+$ put
\begin{equation}\label{set1}\begin{split}
  \hat e_\ga&:=e_\ga\ , \qquad \hat e_{-\ga}=-k_\ga e_{-\ga}\,,\\
  \check e_{-\ga}&:=e_{-\ga}\ , \qquad \check e_{\ga}=- e_{\ga}k_\ga^{-1}\,.
\end{split}\end{equation}

\begin{proposition} For any $\a,\b\in \widehat\Gam$, such that $\a\prec_c\b$,
\begin{align}
[\hat{e}_{\a},\hat{e}_{\b}]_{q^{-1}}&=\sum_{\{\ga_j\},\{n_j\}}
 {C'}_{\{n_j\}}^{\{\ga_j\}}(q)\
\hat{e}^{n_1}_{\ga_1}\hat{e}^{n_2}_{\ga_2}\cdots \hat{e}^{n_m}_{\ga_m},
&& \text{if} \quad \b\in\widehat{\Sigma}_+ \label{use-pro1},\\
[\check{e}_{\a},\check{e}_{\b}]_{q^{-1}}&=\sum_{\{\ga_j\},\{n_j\}}
 {C''}_{\{n_j\}}^{\{\ga_j\}}(q)\
\check{e}^{n_1}_{\ga_1}\check{e}^{n_2}_{\ga_2}\cdots
 \check{e}^{n_m}_{\ga_m},
&& \text{if} \quad \b\in -\widehat{\Sigma}_+ \label{use-pro2},
\end{align}
where  the sums  in \r{use-pro1} and  \r{use-pro2} are
 over all $\ga_1,\ga_2,\ldots,\ga_m$,
$n_1,n_2,\ldots,n_m$ such that
 $\a\ord_c\ga_1\ord_c\ga_2\ord_c\cdots\ord_c\ga_m\ord_c\b$ (meaning
 $\a\prec_c \gamma_1$, $\gamma_1\prec_c \gamma_2$, etc.) and
$\sum_j n_j\ga_j=\a+\b$; $ {C'}_{\{n_j\}}^{\{\ga_j\}}(q)$ and
$ {C''}_{\{n_j\}}^{\{\ga_j\}}(q)$ are Laurent polynomials of $q$.
\end{proposition}

{\em Proof.} See the Appendix. \hfill \qed \medskip

\begin{figure}[t]
\setlength{\unitlength}{0.0007in}
\begin{picture}(9104,3729)(0,-10)
\put(1770,3040){ {$\bullet$}}
\put(2550,2300){ {$\bullet$}}
\put(2620,2100){ {$\bullet$}}
\put(2660,1900){ {$\bullet$}}
\put(2470,1040){ {$\bullet$}}
\put(960,462)  { {$\bullet$}}
\put(2280,2710){ {$\bullet$}}
\put(480,700)  { {$\bullet$}}
\put(-5,1500)  { {$\bullet$}}
\put(60,1300)  { {$\bullet$}}
\put(-10,1720) { {$\bullet$}}
\put(240,2577) { {$\bullet$}}
\put(1800,3200){ {$\gamma_1$}}
\put(2380,2800){ {$\gamma_2$}}
\put(600,282)  { {$-\gamma_1$}}
\put(15,642)   { {$-\gamma_2$}}
\put(2550,850) { {$\gamma_N$}}
\put(-50,2750) { {$-\gamma_N$}}
\put(880,3050) { {$\bullet$}}
\put(670,2950) { {$\bullet$}}
\put(1500,420) { {$\bullet$}}
\put(1140,3100){ {$\bullet$}}
\put(1750,490) { {$\bullet$}}
\put(1970,590) { {$\bullet$}}
\put(6450,2127){ {$\Sigma_{+,e}$}}
\put(6500,1497){ {$\Sigma_{+,F}$}}
\put(5850,1362){ {$\Sigma_{-,f}$}}
\put(5700,1902){ {$\Sigma_{-,E}$}}
\put(6720,2700){ {\small$\{\gamma+m_1\delta\}$}}
\put(7300,1587){ {\small$\{m_2\delta\}$}}
\put(6900,900) { {\small$\{(m_3+1)\delta-\gamma\}$}}
\put(5100,597) { {\small$\{-\gamma-m_4\delta\}$}}
\put(4600,1857){ {\small$\{-m_5\delta\}$}}
\put(4100,2600){ {\small$\{-(m_6+1)\delta+\gamma\}$}}
\put(6100,3500){ {\small1}}
\put(8020,2400){ {\small2}}
\put(6700,2500){ {$\bullet$}}
\put(6500,2570){ {$\bullet$}}
\put(6880,2370){ {$\bullet$}}
\put(5870,1000){ {$\bullet$}}
\put(5700,1120){ {$\bullet$}}
\put(6070,940) { {$\bullet$}}
\put(7110,1750){ {$\bullet$}}
\put(7080,1550){ {$\bullet$}}
\put(6700,1030){ {$\bullet$}}
\put(6900,1170){ {$\bullet$}}
\put(5420,1700){ {$\bullet$}}
\put(5420,1900){ {$\bullet$}}
\put(5650,2350){ {$\bullet$}}
\put(5850,2500){ {$\bullet$}}
\put(1800,2800){\vector(1,-1){500}}
\put(1200,700){\vector(-1,1){500}}
\put(4500,1300){\vector(4,1){3670}}
\put(6700,300){\vector(-1,4){800}}
\end{picture}
\caption{\footnotesize{
The left figure shows the circle ordering of the system
$\widehat\Gam$ of $\Uqg$. It is such that
$\gamma_1+m_1\delta \ord_c m_2\delta \ord_c (m_3+1)\delta-\gamma_2
\ord_c -\gamma_3 -m_4\delta \ord_c -m_5\delta \ord_c
-(m_6+1)\delta+\gamma_4 \ord_c \gamma_5 + m_1\delta$,
where $\gamma_i\in\Sigma_+$ and $m_i>0$. In the right figure,
line 1 shows the decomposition
$\widehat\Gam=\widehat\Sigma_+\sqcup (-\widehat\Sigma_+)$,
related to the subalgebras
$U_q(\mathfrak{n}_\pm)$. Line 2 shows the decomposition
$\Gam=\Gam_E\sqcup \Gam_F$, related to the
 currents Borel subalgebras $U_E$ and $U_F$ discussed in
 Section~\ref{bor-sub}.}}
\label{fig1}
\end{figure}

For $\alpha,\beta\in\widehat\Sigma$ such that $\alpha\prec_c\beta$,
define $[\alpha,\beta] := \{\gamma\in\widehat\Sigma | \alpha\preceq_c\gamma
\preceq_c\beta\}$.

We then define $U^{[\a,\b]}_q(\widehat{\ggg}) \subset \Uqg$ as the subalgebra
generated by the $\hat{e}_\gamma$, $\gamma\in[\a,\b]$,
if $\b\in\widehat\Sigma_+$, and as the subalgebra
generated by the $\check{e}_\gamma$, $\gamma\in[\a,\b]$,
if $\b\in-\widehat\Sigma_+$.

If $\alpha,\beta\in\widehat\Sigma$ are such that $\alpha\prec_c\beta$,
we define the intervals $[\alpha,\beta) := \{\gamma\in\widehat\Sigma |
\alpha\preceq_c\gamma\prec_c \beta \}$ and $(\alpha,\beta]
:= \{\gamma\in\widehat\Sigma | \alpha\prec_c\gamma\preceq_c \beta \}$.
We then define the subalgebras $U^{[\a,\b]}_q(\widehat{\ggg})$,
$U^{[\a,\b)}_q(\widehat{\ggg})$ and $U^{(\a,\b]}_q(\widehat{\ggg})
 \subset \Uqg$
as above (using $\hat e_\gamma$ or $\check e_\gamma$ depending on whether
$\beta\in\widehat\Sigma_+$ or $\beta\in-\widehat\Sigma_+$).

Relations \rf{use-pro1}, \rf{use-pro2} imply that the algebra
$U^{[\alpha,\beta]}_q(\widehat{\ggg})$ admits a
Poincar\'e\--Birkhoff\--Witt (PBW) basis, formed by the ordered monomials
$\hat{e}^{n_1}_{\ga_1}\hat{e}^{n_2}_{\ga_2}\cdots \hat{e}^{n_m}_{\ga_m}$,
where $\alpha = \gamma_1\prec_c\gamma_2\prec_c\ldots \prec_c
\gamma_m = \beta$, if $\beta\in\widehat\Sigma_+$, and by the monomials
$\check{e}^{n_1}_{\ga_1}\check{e}^{n_2}_{\ga_2}\cdots \check{e}^{n_m}_{\ga_m}$,
where $\a=\gamma_1\prec_c\gamma_2\prec_c\ldots \prec_c \gamma_m=\b$, if
$\beta\in-\widehat\Sigma_+$. Moreover, if $\a,\b,\gamma\in \widehat{\Sigma}$
are such that $\a\prec_c\b$, $\b\prec_c\gamma$ and $\a\prec_c\gamma$,
and either $\b,\gamma\in\widehat\Sigma_+$ or
$\beta,\gamma\in-\widehat\Sigma_+$, then the product of $\Uqg$
induces vector space isomorphisms
\begin{equation}\label{Gam0}\begin{split}
U^{[\a,\b)}_q(\widehat{\ggg})\ot U^{[\b,\gamma]}_q(\widehat{\ggg})\simeq
U^{[\a,\gamma]}_q(\widehat{\ggg}),\qquad
U^{[\a,\b]}_q(\widehat{\ggg})\ot U^{(\b,\gamma]}_q(\widehat{\ggg})\simeq
U^{[\a,\gamma]}_q(\widehat{\ggg}),\\
U^{[\a,\b)}_q(\widehat{\ggg})
\otimes U^{[\b,\gamma]}_q(\widehat{\ggg})
\simeq
U^{[\a,\gamma]}_q(\widehat{\ggg}),\qquad
U^{[\a,\b]}_q(\widehat{\ggg}) \otimes
U^{(\b,\gamma]}_q(\widehat{\ggg})
\simeq U^{[\a,\gamma]}_q(\widehat{\ggg}).
\end{split}
\end{equation}

Then $U_q(\nn_+)$ coincides with
$U_q^{\widehat\Sigma_+}(\widehat{\ggg})=U_q^{[\a_{i_1},\delta-\theta]}
(\widehat{\ggg})$,
where $\a_{i_1}$ is the first root (necessarily simple) of the normal
ordering $\prec$ of $\widehat\Sigma_+$. The algebra $U_q(\nn_-)$
coincides with $U_q^{-\widehat\Sigma_+}(\widehat{\ggg})=
U_q^{[-\a_{i_1},-\delta+\theta]}(\widehat{\ggg})$.

The segment $\widehat\Sigma_+=[\a_{i_1},\delta-\theta]$ is a disjoint union
\begin{equation*}
\widehat\Sigma_+ =\Gam_{+,e}\sqcup \Gam_{+,F}\,.
\end{equation*}
Here
\begin{equation}\label{Gam2}
\Gam_{+,e} = \{\gamma + m\delta | \gamma\in \Sigma_+, m\geq 0\},
\quad
\Gam_{+,F} = \{m\delta | m>0\} \sqcup \{-\gamma + m\delta |
\gamma\in\Sigma_+, m>0\}.
\end{equation}
Analogously, $-\widehat\Sigma_+=[-\a_{i_1},-\delta+\theta]$
is a disjoint union
\begin{equation*}
-\widehat\Sigma_+=\Sigma_{-,E}\sqcup \Sigma_{-,f}\,,
\end{equation*}
where
\begin{equation}\label{Gam4}
\Gam_{-,E}= - \Sigma_{+,F}, \quad
\Sigma_{-,f} = -\Sigma_{+,e}.
\end{equation}
The segments \rf{Gam2} and \rf{Gam4} can be united in different ways, composing
segments $\Gam_E$ and $\Gam_F$:
\begin{equation}\label{Gam5}
\begin{split}
\Gam_E&=\Gam_{+,e}\sqcup\Gam_{-,E}=\{n\d, \gamma+m\delta |\,
 \gamma\in \Sigma_+, n< 0,m\in\ZZ\},\\
\Gam_F&=\Gam_{+,F}\sqcup\Gam_{-,f}=\{n\d, -\gamma+m\delta |\,
 \gamma\in \Sigma_+, n> 0,m\in\ZZ\}.
\end{split}
\end{equation}
The subalgebras related to the segments \rf{Gam2}, \rf{Gam4} and
\rf{Gam5} play a crucial role in our further constructions.

\subsection{The `currents' presentation of $\Uqg$}\label{section2.3}

In this presentation (\cite{D88}), $\Uqg$ is generated by the central
elements $\C^{\pm 1}$, the grading elements $q^{\pm d}$, and by the elements
$\ein{\a}{n}$, $\fin{\a}{n}$ (where $\a\in\Pi$, $n\in\ZZ$) and
$k_\a^{\pm 1}$, $\han{\a}{n}$ (where $n\in\ZZ\,\backslash \{0\}$,
$\a\in\Pi$).
These elements are gathered into generating functions
\begin{eqnarray*}
&\ds e_\a(z)\,=\,\sum_{n\in\z}\ein{\a}{n}z^{-n}\ ,
\quad f_\a(z)\,=\,\sum_{n\in\z}\fin{\a}{n}z^{-n}\ , \\
&\ds\psi^\pm_\a(z)\,=\sum_{n>0}\psi_\a^{\pm}[n]z^{\mp n}
=k_{\a}^{\pm1}
\,
\exp\sk{\pm(q_\a-q_\a^{-1})\sum_{n>0}\hin{\a}{\pm n}z^{\mp n}}\ ,
\end{eqnarray*}
such that $q^d q^{-d} = q^{-d} q^d = CC^{-1} = C^{-1}C = 1$ and
$q^da(z)q^{-d}=a(q^{-1}z)$ for any of these generating
functions, which we will call currents. Currents are labeled by the
simple roots $\a\in\Pi$ of $\mathfrak{g}$.

The defining relations of $\Uqg$ in the `currents' presentation are
($\alpha,\beta\in\Pi$):
\begin{equation}
\begin{split}
(z-q^{({\a},{\b})}w)e_{\a}(z)e_{\b}(w)&=
 e_{\b}(w)e_{\a}(z)(q^{({\a},{\b})}z-w)\ ,\\
(z-q^{-({\a},{\b})}w)f_{\a}(z)f_{\b}(w)&=
 f_{\b}(w)f_{\a}(z)(q^{-({\a},{\b})}z-w)\ , \\
\label{1}
\ds
 \psi_{\a}^\pm(z)e_{\b}(w)\left(\psi_{\a}^\pm(z)\right)^{-1}&=\ds
\frac{q^{({\a},{\b})}\C^{\pm 1}z-w}
{\C^{\pm 1}z- q^{({\a},{\b})}w}e_{\b}(w)\ ,\\
\ds
\psi_{\a}^\pm(z)f_{\b}(w)\left(\psi_{\a}^\pm(z)\right)^{-1}&=\ds
\frac{q^{-({\a},{\b})}\C^{\mp 1}z-w}
{\C^{\mp 1}z- q^{-({\a},{\b})}w}f_{\b}(w)\ ,
 \\
\psi_{\a}^\pm(z)\psi_{\b}^\pm(w)&=
 \psi_{\b}^\pm(w)\psi_{\a}^\pm(z)\ ,
\\
\frac{q^{({\a},{\b})}z-\C^{2}w}
{z- q^{({\a},{\b})}\C^{2}w}\psi_\a^+(z)\psi_\b^-(w)&=
\frac{q^{({\a},{\b})}\C^{2}z-w}
{\C^{2}z- q^{({\a},{\b})}w}\psi_\b^-(w)\psi_\a^+(z)\ ,
\\
[e_{\a}(z),f_{\b}(w)]=
 \frac{\delta_{\a,\b}}{q_{\a}-q^{-1}_{\a}}\left(\delta(z/\C^2w)\right.&\left.
\psi^+_{\a}(\C^{-1}z)-\delta(\C^2z/w)\psi^-_{\a}(\C^{-1}w)\right)\,,
\end{split}
\end{equation}
\begin{equation*}
\sum_{r=0}^{n_{ij}}(-1)^r\left[\gfr{n_{ij}}{r}\right]_{q_{i}}
{\rm Sym}_{z_1,...,z_{n_{ij}}}e_{\pm\al_i}(z_1)\cdots e_{\pm\al_i}(z_r)
e_{\pm\al_j}(w)e_{\pm\al_i}(z_{r+1})\cdots e_{\pm\al_i}(z_{n_{ij}})=0.
\end{equation*}
Here $\a_i\not=\a_j$, $n_{ij}=1-a_{i,j}$, $\delta(z)=\sum_{k\in\ZZ}z^k$.

We now describe the isomorphism of the two realizations (\cite{D88,CP}).
Suppose that the root vector  $e_\theta\in\ggg$, corresponding to
the longest root $\theta\in \Sigma_+$,
is presented as a multiple commutator
$e_\theta= \lambda[e_{\a_{i_1}},[e_{\a_{i_2}},
\cdots [e_{\a_{i_n}},e_{\a_{j}}]\cdots]]$
for some $\lambda\in\CC$. The isomorphism is given by
the assignment
\begin{equation}
\begin{split} \label{direct:isom}
& k_{\a_i}^{\pm 1}\mapsto k_{\a_{i}}^{\pm 1}\,,\qquad
e_{\a_i}\mapsto \ein{\a_i}{0}\,,\qquad
e_{-\a_i}\mapsto \fin{\a_i}{0}\,,\qquad i=1,\ldots,r,\\
&k_{\a_0}^{\pm 1}\to \C^{\pm 2}k_\theta^{\mp 1}=\C^{\pm 2}\prod_{i=1}^r
k_{\a_i}^{\mp n_i},
\qquad \quad k_\delta^{\pm 1/2} \mapsto C^{\pm 1},
\qquad \quad q^{\pm d}\mapsto q^{\pm d},
\\
& e_{\a_0}\mapsto \mu S^-_{i_1}S^-_{i_2}\cdots S^-_{i_n}(\fin{\a_j}{1})\,,\qquad
e_{-\a_0}\mapsto \lambda S^+_{i_1}S^+_{i_2}\cdots S^+_{i_n}(\ein{\a_j}{-1})
\end{split}
\end{equation}
Here $S_{i}^\pm\in \on{End}(U_q(\widehat{\ggg}))$ are the following
operators of adjoint action (with respect to the coproducts $\Delta^{std}$
and $(\Delta^{std})^{2,1}$, see \rf{copr}):
\begin{equation}\label{ad}
S_i^+(x)= \ein{\a_i}{0}x-k_{\a_i}xk_{\a_i}^{-1}\ein{\a_i}{0},\qquad
S_i^-(x)= x\fin{\a_i}{0}-\fin{\a_i}{0}k_{\a_i}{}xk_{\a_i}^{-1}\ ,
\end{equation}
and the constant $\mu$ is determined by the condition that relation
\rf{Q1} for $i=j=0$ remains valid in the image.

We now describe the inverse isomorphism.
Let $\pi:\{\a_1,\ldots,\a_r\}\mapsto \{0,1\}$ be the map such that
$\pi(\a_i)\neq\pi(\a_j)$ for $(\a_i,\a_j)\neq0$ and $\pi(\alpha_1)=0$.
\begin{proposition}\label{KhT-th} (see \cite{KhT,Be,Dam})
The inverse of isomorphism (\ref{direct:isom}) is such that
\begin{equation}\label{UEgen}
\begin{split}
e_{\a_i}[n]& \mapsto (-1)^{n\pi(\a_i)}\hat{e}_{\a_i+n\d}=\left\{
\begin{array}{l}
\ds (-1)^{n\pi(\a_i)} e_{\a_i+n\d}\,,\quad n\geq 0,\\
\ds -(-1)^{n\pi(\a_i)}k_{-\a_i-n\d}e_{\a_i+n\d}\,,\quad n< 0,
\end{array}
\right. \\
f_{\a_i}[n]& \mapsto (-1)^{n\pi(\a_i)}\check{e}_{-\a_i+n\d}=\left\{
\begin{array}{l}
\ds (-1)^{n\pi(\a_i)} e_{-\a_i+n\d}\,,\quad n\leq 0,\\
\ds - (-1)^{n\pi(\a_i)} e_{-\a_i+n\d}\, k^{}_{\a_i-n\d}\,,\quad n>0 ,
\end{array}
\right.
\end{split}
\end{equation}
\begin{equation}\label{UPp}
\begin{split}
\psi^+_{\a_i}[0]& \mapsto k_{\a_i}\,,\qquad
\psi^-_{\a_i}[0] \mapsto k^{-1}_{\a_i}\,,\qquad
C^{\pm 1}\mapsto k_\d^{\pm 1/2},
\\
\psi^+_{\a_i}[n]&\mapsto (q-q^{-1})(-1)^{n\pi(\a_i)}
k_{\a_i}k_\d^{-\frac{n}{2}}e^{\prime(i)}_{n\d}\,,\quad n> 0,\\
\psi^-_{\a_i}[n]&\mapsto -(q-q^{-1})(-1)^{n\pi(\a_i)}
k^{-1}_{\a_i}k_\d^{-\frac{n}{2}} e^{\prime(i)}_{n\d}\,,\quad n< 0.
\end{split}
\end{equation}
\end{proposition}
The relation between the imaginary root generators $e^{\prime(i)}_{k\d}$ and
$e^{(i)}_{k\d}$ is given by formulas \r{im-rel}.

Note that the root vectors $e_{\pm\a_i+n\d}$, $n\in\ZZ$, as well
as the imaginary root generators do not depend on the choice of a normal
ordering of $\widehat\Sigma_+$, satisfying the condition \rf{ord1} (\cite{Dam},
Proposition 11). So the identification of Proposition \ref{KhT-th}
 does not depend on such a normal ordering.

The `currents' bialgebra structure $\Delta^{(D)}$ on $\Uqtri$ is given by:
\begin{align}
\notag
\Delta^{(D)}(q^{\pm d})&=q^{\pm d}\ot q^{\pm d}, \quad
\Delta^{(D)}(\C^{\pm 1})=\C^{\pm 1}\ot \C^{\pm 1},\\
\Delta^{(D)}\psi_\a^\pm(z)&=\ \psi_\a^\pm(\C_2^{\pm 1}z)\otimes
\psi_\a^\pm(\C_1^{\mp 1}z)\, , \label{compsi}\\
\Delta^{(D)}\xp_\a(z)&=\ \xp_\a(z)\otimes 1+\psi_\a^-(\C_1 z)\otimes \xp_\a(\C_1^2z)\,
\label{come} ,\\
\Delta^{(D)}\xm_\a(z)&=\ 1\otimes \xm_\a(z)+\xm_\a(\C_2^2z)\otimes
\psi_\a^+(\C_2z)\, ,\label{comf}
\end{align}
where $C_1=C\otimes1$ and $C_2=1\otimes C$.
The counit map is given by:
$$
\coun(e_\a(z))=\coun(f_\a(z))=0\,,\quad
\coun(\psi_\a^\pm(z))=\coun(q^{\pm d}) = \coun(C^{\pm 1}) =1.
$$

The principal degree on $\Uqtri$ is defined by $|e_i[n]| = n\nu +1$,
$|k_i^{\pm 1}| = |C| = 0$, $|f_i[n]|=n \nu -1$, $|h_\alpha[n]|=
n\nu$, where $\nu = \sum_{i=0}^r n_i$ (recall that $\delta = \sum_{i=0}^r
n_i \alpha_i$); then $\Delta^{(D)} : \Uqtri \to
\Uqtri^{\otimes_< 2}$ is a topological bialgebra, in the sense of
Section \ref{assum}, (H1). The identification with the situation of
this Section is $D = \Uqtri$, $\Delta_D = \Delta^{(D)}$.

The coproducts $\Delta^{std}$ of Subsection~\ref{ch-gen} and $\Delta^{(D)}$ are
related by a twist which can be described explicitly (see
Proposition \ref{prop:nontriv} and \cite{KhT}).

\subsection{Intersections of Borel subalgebras}\label{bor-sub}

We first describe the Borel subalgebras related to the `currents'
realization of $\Uqg$.

Let $U_F$ be the subalgebra of $\Uqg$ generated by $U_q(\mathfrak{h})$ and
the elements $\fin{\a}{n}$ (where $\a\in\Pi$, $n\in\ZZ$) and
$\hin{\a}{m}$ (where $\a\in \Pi$, $m>0$). In the circular description,
this is the subalgebra of $\Uqg$ associated with the segment
$\Gam_F$ (see \rf{Gam5}). Analogously, we define
$U_E$ as the subalgebra of $\Uqg$ generated by $U_q(\mathfrak{h})$
and the elements $\ein{\a}{n}$ (where $\a\in\Pi$, $n\in\ZZ$) and
$\hin{\a}{m}$ (where $\a\in\Pi$, $m<0$). It corresponds to the subalgebra
of $\Uqg$ associated to the segment $\Gam_E$. Both $U_F$ and $U_E$
are Hopf subalgebras of ${\Uqg}$ with respect to the
coproduct $\Delta^{(D)}$. We call them the {\it currents Borel subalgebras}.
We also define $U_f \subset U_F$ as the subalgebra generated by the elements
$\fin{\a}{n}$ (where $\a\in\Pi$, $n\in\ZZ$) and $U_e \subset U_E$ as the
subalgebra generated by the elements $\ein{\a}{n}$
(where $\a\in\Pi$, $n\in\ZZ$).

We have
$$
\Uqtri \simeq U_E \otimes_{U_q({\mathfrak h})} U_F,
$$
where the isomorphism is induced by the product map.
So $D = \Uqtri$, $A := U_E$, $B := U_F$, $\CC[\ZZ^r] := U_q({\mathfrak h})$
satisfy the first part of (H2) in Section \ref{assum}.

We define $U_F^0\subset U_F$, $U_E^0\subset U_E$ as the algebras
with the same generators, except $U_q(\h)$.

The relations \r{come} and \r{comf} imply that the subalgebra $U_f$
is a right coideal of $U_F$ with respect to $\Delta^{(D)}$, and the
subalgebra $U_e$ is a left coideal of $U_E$ with respect to
$\Delta^{(D)}$:
\begin{equation} \label{coprod:D}
\Delta^{(D)}(U_f)\subset U_f \ot U_F, \qquad
\Delta^{(D)}(U_e)\subset U_E \ot U_e\ .
\end{equation}

There is a unique bialgebra pairing $\langle - , -\rangle : U_E \otimes U_F
\to \CC$, expressed by
$$
\langle e_\alpha(z), f_\beta(w) \rangle = {\delta_{\alpha,\beta}
\delta(z/w)\over {q_\alpha^{-1} - q_\alpha}}, \quad
\langle \psi_\alpha^-(z),\psi_\beta^+(w) \rangle =
{{q^{(\alpha,\beta)} - z/w}\over{1-q^{(\alpha,\beta)}z/w}},
\quad \langle c,d \rangle=1
$$
in terms of generating series; all other pairings
between generators are zero.

\begin{proposition}
This pairing satisfies identity (\ref{D-mu}) (see (H2)).
\end{proposition}

{\em Proof.} Let us set $R(a,b) :=
a^{(1)} b^{(1)} \langle a^{(2)},b^{(2)} \rangle
- b^{(2)} a^{(2)} \langle a^{(1)},b^{(1)} \rangle$ for
$a\in U_E$, $b\in U_F$. One checks that $R(a,b) = 0$ if
$a,b$ are generators of $U_E$, $U_F$. Moreover, we have the identities
$$
R(aa',b) = a^{(1)}R(a',b^{(1)}) \langle a^{(2)},b^{(2)}\rangle
+ R(a,b^{(2)}) a^{\prime(2)} \langle a^{\prime(1)}, b^{(1)} \rangle ,
$$
$$
R(a,bb') = R(a^{(1)},b) b^{\prime(1)} \langle a^{(2)},b^{\prime(2)}\rangle
+ b^{(2)} R(a^{(2)},b') \langle a^{(1)}, b^{(1)} \rangle .
$$
Reasoning by induction on the length of $a$ and $b$ (expressed as
products of generators), these identities imply that $R(a,b) = 0$
for any $a,b$. \hfill \qed \medskip

Therefore $\langle-,-\rangle$ satisfies the hypothesis (H2).
One also checks that it satisfies hypothesis (H6).

We will be interested in intersections of Borel subalgebras
of different types. Denote by $U_F^+$, $U_f^-$, $U_e^+$ and $U_E^-$ the
following intersections of the standard and currents Borel
algebras:
\begin{align}
U_f^-=U_F\cap U_q(\mathfrak{n}_-), \quad
U_F^+ =U_F\cap U_q(\mathfrak{b}_+), \quad
U_F^{0,+} =U_F^0\cap U_q(\mathfrak{n}_+),
\nonumber \\ 
U_e^+=U_E\cap U_q(\mathfrak{n}_+), \quad
U_E^-=U_E\cap U_q(\mathfrak{b}_-), \quad
U_E^{0,-}=U_E^0\cap U_q(\mathfrak{n}_-).
\end{align}
The notation is such that the upper sign says in which standard Borel
subalgebra $U_q(\mathfrak{b}_\pm)$ the given algebra is contained and
the lower letter says in which currents Borel subalgebra ($U_F$ or $U_E$)
it is contained. This letter is capital if the subalgebra contains
imaginary root generators $\hin{i}{n}$ and the Cartan subalgebra
$U_q(\mathfrak{h})$ and is small otherwise. In the notation of Section
\ref{sectionCW}, the subalgebras $U_f^-$ and $U_F^+$ correspond to the
segments $\Gam_{-,f}$ and $\Gam_{+,F}$; the subalgebras $U_E^-$ and $U_e^+$
correspond to the segments $\Gam_{-,E}$ and $\Gam_{+,e}$.

\begin{proposition}\label{subalgebras}
The product in $U_E$ sets up an isomorphism of vector
spaces $U_E \simeq U_e^+\ot U_E^-$\,.
The product in $U_F$ sets up an isomorphism of vector
spaces $U_F \simeq U_f^-\ot U_F^+$\,.
\end{proposition}
\noindent
{\it Proof}. This follows from (\ref{Gam0}).
\hfill$\Box$ \medskip

We will set $A_1 := U_e^+$, $A_2 := U_E^-$, $B_1 := U_F^+$, $B_2 :=
U_f^-$. Then the hypothesis (H3) is satisfied.
We will set $B'_1 := U_F^{0,+}$, $A'_2 := U_E^{0,-}$, so that
(H5) is satisfied.

The next proposition describes a family of generators of the intersections
of Borel subalgebras.
\begin{proposition}${}$\label{subalgebras2}
\begin{itemize}
\item[(i)]
The algebra  $U_e^+$  (resp., $U_f^-$) is generated by the elements
$\ein{i}{n}$ (resp., $\fin{i}{-n}$), where $i\in \{1,\ldots,r\}$, $n\geq0$.

\item[(ii)]
The algebra $U_E^-$ (resp., $U_F^+$) is generated by
$U_q({\mathfrak h})$, the elements
$\ein{i}{n}$, $\hin{i}{m}$, where $n,m<0$, $i\in \{1,\ldots,r\}$, and
by the root vectors $e_{\gamma-\delta}$, $\gamma\in\Sigma_+$
(resp., $U_q({\mathfrak h})$,
the elements $\fin{i}{n}$, $\hin{i}{m}$,
where  $n, m>0$, $i\in \{1,\ldots,r\}$,
and the root vectors $e_{\delta-\gamma}$, $\gamma\in\Sigma_+$).
\end{itemize}
\end{proposition}

\noindent
{\it Proof.} The PBW result shows that a basis for $U_e^+$ is given by the
ordered monomials in the $\hat e_\gamma$, $\gamma\in \Sigma_{+,e}$.
These generators are expressed via those listed in the Proposition,
by means of successive applications of relations \rf{use-pro1} and
\rf{use-pro2}. The proof is the same in the other cases.
\hfill$\Box$ \medskip

Note that the generators listed in the Proposition do not depend
on a choice of the normal ordering, satisfying \rf{ord1}. This was
already stated in Section \ref{sectionCW} for all of them, except
for the root vectors $e_{\pm(\delta-\gamma)}$, $\gamma\in\Sigma_+$. They are
constructed as successive $q$-commutators of the type
$[e_{\pm\gamma'},e_{\pm(\delta-\gamma'')}]_{q^{\pm 1}}$. By induction of
the height of $\gamma$, one proves that these $q$-commutators define the
same root vectors.

\subsection{Relation between $\Delta^{std}$ and $\Delta^{(D)}$}
\label{sect:rel}

\begin{proposition} \label{prop:nontriv}
We have $\Delta^{(D)}(U_q({\mathfrak b}_-)) \subset U_q({\mathfrak b}_-)
\otimes U_q({\mathfrak g})$ and $\Delta^{(D)}(U_q({\mathfrak b}_+))
\subset U_q({\mathfrak g}) \otimes U_q({\mathfrak b}_+)$.
\end{proposition}

{\em Proof.} We prove the first statement.
It suffices to show that $\Delta^{(D)}(e_{-\alpha_i}) \in
U_q({\mathfrak b}_-) \otimes U_q({\mathfrak g})$, for $i = 0,...,r$.
When $i\in \{1,...,r\}$, $\Delta^{(D)}(e_{-\alpha_i}) = \Delta^{(D)}
(f_{\alpha_i}[0]) = 1\otimes \fin{\a_i}{0}+\sum_{k\geq0}
\fin{\a_i}{-k}\otimes \psi^+_{\a_i}[k]C^{-2n+k} \in
U_q({\mathfrak b}_-) \otimes U_q({\mathfrak g})$.

When $i=0$, $\Delta^{(D)}(e_{-\alpha_0})=\lambda\Delta^{(D)}(S_{i_1}^+...
S_{i_n}^+(e_{\alpha_j}[-1]))$.
Let us prove by descending induction on $k$ that
$\Delta^{(D)}(S_{i_k}^+...S_{i_n}^+(e_{\alpha_j}[-1])) \in
U_E^- \otimes U_q({\mathfrak g})$. This is true for
$k=n+1$ since $\Delta^{(D)}(e_{\alpha_j}[-1]) = e_{\alpha_j}[-1] \otimes 1
+ \sum_{s\geq 0} \psi^-_{\alpha_j}[-s] C^{2-s}\otimes e_{\alpha_j}[s-1]$.
If now $x = S_{i_{k+1}}^+...S_{i_n}^+(e_{\alpha_j}[-1])$ is such that
$\Delta^{(D)}(x)\in U_E^- \otimes U_q({\mathfrak g})$, we have
$\Delta^{(D)}(S_i^+(x)) = \Delta^{(D)}(e_i[0]x - q^{(\alpha_i,|x|)}xe_i[0])
= [\Delta^{(D)}(e_{\alpha_i}[0]),x]_q = [e_{\alpha_i}[0] \otimes 1
+ \sum_{s\geq 0} \psi^-_{\alpha_i}[-s]C^{-s} \otimes e_{\alpha_i}[s],
\Delta^{(D)}(x)]_q \equiv [e_{\alpha_i}[0] \otimes 1,\Delta^{(D)}(x)]_q$
modulo $U_E^- \otimes U_q({\mathfrak g})$ (here $|x|$ is the degree of $x$).
The result now follows from the induction assumption $\Delta^{(D)}(x)
\in U_E^- \otimes U_q({\mathfrak g})$ and (\ref{use-pro1}).

The proof of the second statement is similar.
\hfill \qed \medskip

\begin{corollary}
The intersections of Borel subalgebras have the following
coideal properties:
  \begin{align}
  \Delta^{(D)}(U_F^+)&\subset U_F\ot_< U_F^+, &
  \Delta^{(D)}(U_f^-)&\subset U_f^-\ot_< U_F, \label{coid1}\\
  \Delta^{(D)}(U_E^-)&\subset U_E^-\ot_< U_E, &
  \Delta^{(D)}(U_e^+)&\subset U_E \ot_< U^+_e, \label{coid2}
  \end{align}
\end{corollary}

{\em Proof.} This follows directly from (\ref{coprod:D}) and
Proposition \ref{prop:nontriv}.
\hfill \qed \medskip

This means that $(D,\Delta^{(D)})$, $A_i$, $B_i$ satisfy the hypothesis
(H4).

We can therefore apply Theorem \ref{thm:twists} to
to $(D,\Delta_D) = (U_q(\widehat{{\mathfrak g}}), \Delta^{(D)})$,
$A = U_E$, $A_1 = U_e^+$, $A'_2 = U_E^{0,-}$,
$B = U_F$, $B'_1 = U_F^{0,+}$, $B_2 = U_f^-$, $\CC[\ZZ^r] = U_q(\h)$.


We denote by ${\cal R}_1\in U_e^+ \otimes U_f^-$ and
${\cal R}_2\in U_E^- \otimes U_F^+$ the analogues of $R_1$, $R_2$.
Then ${\cal R}_2$ is a cocycle for $(U_q(\widehat{\mathfrak g}),\Delta^{(D)})$,
so we get a Hopf algebra structure on $U_q(\widehat{\mathfrak g})$,
given by $\Delta^{tw}(x) = {\cal R}_2 \Delta^{(D)}(x) {\cal R}_2^{-1}$
($\Delta^{tw}$ is the analogue of $\Delta$ of section \ref{assum}),
quasitriangular with $R$-matrix ${\cal R}^{tw} = {\cal R}_2^{2,1}{\cal R}_1$.

A proof of the following result (based on the braid group action)
was first given in \cite{KhT}.

\begin{proposition} 
$\Delta^{tw} = (\Delta^{std})^{2,1}$, therefore we have $\Delta^{(D)}(x)
= {\cal R}_2^{-1} (\Delta^{std})^{2,1}(x) {\cal R}_2$. We also have
\begin{equation}\label{twist:rel}
\Delta^{(D)}(x) = {\cal R}_1^{2,1} \Delta^{std}(x)
({\cal R}_1^{2,1})^{-1}.
\end{equation}
\end{proposition}

{\em Proof.} $\Delta^{tw}$ defines a Hopf structure on
$U_q(\widehat{\mathfrak g})$, for which $U_q({\mathfrak b}_+) = A_1B_1$
and $U_q({\mathfrak b}_-) = A_2B_2$ are Hopf subalgebras.
Since ${\cal R}_2$ is invariant under $U_q({\mathfrak h})$, we have
$\Delta^{tw}(k^{\pm 1}_{\alpha_i}) = (k_{\alpha_i}^{\pm 1})^{\otimes 2}$ and
$\Delta^{tw}(k^{\pm 1/2}_{\delta}) = (k_{\delta}^{\pm 1/2})^{\otimes 2}$.

Since $\Delta^{(D)}$ and ${\cal R}_2$ are homogeneous for the
principal degree, the same holds for $\Delta^{tw}$. Moreover,
${\cal R}_2 = {\cal R}_{{\mathfrak h}}(1 + ($negative principal degree$)
\otimes ($positive principal degree$))$, where ${\cal R}_{{\mathfrak h}} =
q^{C_{{\mathfrak h}}}$ and $C_{{\mathfrak h}}\in {\mathfrak h}^{\otimes 2}$
is the Casimir element of the Cartan algebra. Therefore
$\Delta^{tw}(e_{\alpha_i}) = e_{\alpha_i} \otimes k_{\alpha_i}
+ ($degree 0$) \otimes ($degree 1$)$, and $\Delta^{tw}(e_{-\alpha_i}) =
1 \otimes e_{-\alpha_i} + ($degree $-1) \otimes ($degree 0$)$.

On the other hand, $\Delta^{tw}(x) = {\cal R}_2\Delta^{(D)}(x) {\cal R}_2^{-1}
= {\cal R}_1^{-1}(\Delta^{(D)})^{2,1} {\cal R}_1$, and since
${\cal R}_1 = 1 + ($positive principal degree$) \otimes ($negative principal
degree$)$, we find
$\Delta^{tw}(e_{\alpha_i}) = 1 \otimes e_{\alpha_i} + ($degree 1$) \otimes
($degree 0$)$, and $\Delta^{tw}(e_{-\alpha_i}) =
e_{-\alpha_i} \otimes k_{-\alpha_i}+ ($degree 0$) \otimes ($degree $-1)$.
Combining these results, we get $\Delta^{tw} = (\Delta^{std})^{2,1}$

Then $\Delta^{std}(x) = {\cal R}_2^{2,1} (\Delta^{(D)})^{2,1}(x)
({\cal R}_2^{2,1})^{-1} = ({\cal R}_1^{2,1})^{-1} \Delta^{(D)}(x)
{\cal R}_1^{2,1}$, which implies (\ref{twist:rel}).
\hfill \qed \medskip

\medskip

\begin{proposition} \label{coidealprop}
The intersections of Borel subalgebras have the following
coideal properties:
  \begin{align}
 \Delta^{std}(U_F^+)&\subset U_F^+\ot U_q(\mathfrak{b}_+), &
  \Delta^{std}(U_f^-)&\subset U_f^-\ot U_q(\mathfrak{b}_-), \label{coid3}\\
  \Delta^{std}(U_E^-)&\subset U_q(\mathfrak{b}_-)\ot U_E^-, &
  \Delta^{std}(U_e^+)&\subset U_q(\mathfrak{b}_+)\ot U^+_e. \label{coid4}
  \end{align}
\end{proposition}

This is a translation of Proposition \ref{prop:coid:gal}.

We will denote by $P^\pm$ the projection operators of the Borel
subalgebra $U_F$, corresponding to the decomposition $U_F=U_f^-U_F^+$.
So for any $f_+\in U_F^+$ and any $f_-\in U_f^-$,
\begin{align}
 \label{Pdef}
 P^+(f_-f_+)&=\varepsilon(f_-)f_+,&
 P^-(f_-f_+)&=f_-\varepsilon(f_+)\ .
\end{align}
The operator $P^+$ will also be denoted by $P$ (without index).

\setcounter{equation}{0}
\section{Universal weight functions}

\subsection{The definition}\label{sectiondef}

Denote by $\left(U_e^+\right)^\varepsilon$ the augmentation
ideal of $U_e^+$, i.e., $\left(U_e^+\right)^\varepsilon= U_e\cap
\text{Ker}(\varepsilon)$. Let $J$ be the left ideal of $U_q(\bb_+)$
generated by $\left(U_e^+\right)^\varepsilon$:
\begin{equation}\label{ideal}
J=U_q(\bb_+)\left(U_e^+\right)^\varepsilon=
\sum_{\alpha\in\Pi,n\geq 0} U_q(\bb_+) \ein{\a}{n}.
\end{equation}

Recall that the $\psi_{\alpha_i}^+[0]^{\pm 1}$, $\psi_{\alpha_i}^+[n]$
and $C^{\pm 1}$ (where $i\in \{1,...,r\}$, $n>0$) commute in
$U_q(\widehat{\ggg})$.

\begin{proposition}
\label{propJ}
$J$ is a  coideal of $U_q(\bb_+)$, i.e.,
$\Delta^{std}(J)\subset J\ot U_q(\bb_+) + U_q(\bb_+)\ot J$.
The space $J$ is also stable under right multiplication by the
$\psi_{\alpha_i}^+[0]^{\pm 1}$, $\psi_{\alpha_i}^+[n]$
(where $i\in \{1,..,r\}$, $n>0$), i.e.,
$J \psi_{\alpha_i}^+[0]^{\pm 1} \subset J$,
$J \psi_{\alpha_i}^+[n] \subset J$.

It follows that $U_q({\mathfrak b}_+)/J$ is both a coalgebra and a right
module over $\CC[\psi_{\alpha_i}[0]^{\pm 1},\psi_{\alpha_i}[n],$
$C^{\pm 1}; i\in \{1,...,r\}, n>0]$.
\end{proposition}

{\em Proof.} The second part of (\ref{coid4}) implies that for $n\geq 0$,
$\Delta^{std}(e_\alpha[n]) \in U_q({\mathfrak b}_+) \otimes U_e^+$.
Moreover, $(\on{id} \otimes \varepsilon) \circ \Delta^{std} = \on{id}$
implies that $\Delta^{std}(e_\alpha[n]) - e_\alpha[n]\otimes 1 \in
U_q({\mathfrak b}_+) \otimes \on{Ker}\varepsilon$. Therefore
$\Delta^{std}(e_\alpha[n]) - e_\alpha[n]\otimes 1 \in
U_q({\mathfrak b}_+) \otimes (U_e^+)^\varepsilon$, so
$\Delta(e_\alpha[n]) \in J \otimes U_q({\mathfrak b}_+)
+ U_q({\mathfrak b}_+) \otimes J$. This implies the coideal property of $J$.
The second property follows from the commutation relation of the
currents $\psi_{\alpha_i}^+(w)$ and $e_\alpha(z)$.
\hfill \qed \medskip

Define an {\it ordered $\Pi$-multiset} as a triple $\bar I = (I,\prec,\iota)$,
where $(I,\prec)$ is a finite, totally ordered set, and $\iota : I\to \Pi$
is a map (the `coloring map'). Ordered $\Pi$-multisets form a category,
where a morphism $\bar I\to \bar I' = (I',\prec',\iota')$ is a map
$m : I\to I'$, compatible with the order relations and such that
$\iota' \circ m = m \circ \iota$.

If $\bar I = (I,\prec,\iota)$ is an ordered $\Pi$-multiset, then a
partition $I = I_1\sqcup I_2$ gives rise to ordered $\Pi$-multisets
$\bar I_i = (I_i,\prec_i,\iota_i)$, $i = 1,2$, where $\prec_i$ and
$\iota_i$ are the restrictions of $\prec$ and $\iota$ to $I_i$.

If $\bar I = (I,\prec,\iota)$ is an ordered $\Pi$-multiset, and
$I = \{i_1,...,i_n\}$, with $i_1\prec...\prec i_n$, then we attach
to $\bar I$ an ordered set of variables $(t_i)_{i\in I}
=(t_{i_1},...,t_{i_n})$. The `color' of $t_k$ is $\iota(i_k)\in\Pi$.

For any vector space $V$ denote by $V[[t_1^{\pm 1},...,t_n^{\pm 1}]]$
the vector space of all formal series
\begin{equation*}
\sum_{(k_1,...,k_n)\in\ZZ^n} A_{k_1,...,k_n}t_1^{k_1}\cdots t_n^{k_n},
\qquad A_{k_1,...,k_n}\in V\,,
\end{equation*}
and by $V((t_1))\cdots ((t_n))$ its subspace corresponding to the maps
$(k_1,...,k_n)\mapsto A_{k_1,...,k_n}$, such that
there exists an integer $a_1$ and
integer valued functions $a_2(k_1)$, $a_3(k_1,k_2)$, $\ldots$,
$a_n(k_1,...,k_{n-1})$ (depending on the map $(k_1,...,k_n)\mapsto
A_{k_1,...,k_n}$), such that $A_{k_1,...,k_n}=0$ whenever
$k_1>a_1$, or $k_2>a_2(k_1)$, or $k_3>a_3(k_1,k_2)$, $\ldots\ $, or
$k_n>a_n(k_1,...,k_{n-1})$. A bilinear map $V\otimes W\to Z$ gives rise to
a bilinear map $V((t_1))...((t_n)) \times W((t_1))...((t_n)) \to
Z((t_1))...((t_n))$ using the multiplication of formal series.

A {\it universal weight function} is an assignment $\bar I \mapsto
\on{W}_{\bar I}$, where $\bar I = (I,\prec,\iota)$ is an ordered
$\Pi$-multiset and $W_{\bar I}((t_i)_{i\in I}) \in
(U_q({\mathfrak b}_+)/J)((t_{i_n}^{-1}))...((t_{i_1}^{-1}))$
($I = \{i_1,...,i_n\}$,
with $i_1\prec ... \prec i_n$) is such that
\begin{itemize}
\item[(a)] (functoriality) If $f : \bar I \to \bar J$ is an isomorphism of
ordered $\Pi$-multisets, then
$\on{W}_{\bar I}((t_i)_{i\in I})=\on{W}_{\bar J}((t_{f^{-1}(j)})_{j\in J})$;
\item[(b)] $\on{W}_{\emptyset} = 1$, where $\on{W}_{\emptyset}$ is the series
corresponding to $I$ equal to the empty set;
\item[(c)] The series $\on{W}_{\bar I}(t_{i_1},\ldots,t_{i_n})$ satisfies
the relation
\begin{equation}
\begin{split}
&\Delta^{std} \sk{\on{W}_{\bar I}((t_i)_{i\in I})}=
\sum\limits_{I=I_1\sqcup I_2}
\big( \on{W}_{\bar I_1}((C_2^2t_i)_{i\in I_1})
\otimes
\on{W}_{\bar I_2}((t_i)_{i\in I_2}) \big) \times\\&\qquad
\times
\big( 1 \otimes \prod_{i\in I_1} \psi^+_{\i(i)}(Ct_{i}) \big) \times
\prod\limits_{\substack{k,l | k<l, \\ i_k\in I_1,i_l\in I_2}}
\frac{q^{-(\a_{\i(i_k)},\a_{\i(i_l)})}-t_{i_l}/t_{i_k}}
{1-q^{-(\a_{\i(i_k)},\a_{\i(i_l)})}t_{i_l}/t_{i_k}}\\
\label{univ-dec}
\end{split}
\end{equation}
where $\Delta^{std}$ is the coproduct of $U_q({\mathfrak b}_+)/J$.
Here $C_2 = 1\ot C$,
where $C$ is the central element of $\Uqg$, see Section \ref{section2.3}.
The summation in \r{univ-dec} runs over all possible decompositions of the
set $I$ into two disjoint subsets $I_1$ and $I_2$.
\end{itemize}

\bigskip

\noindent {\bf Remark}. Let $\Uqgp$ be the subquotient of $\Uqg$ with
trivial central element $C=1$ and dropped gradation element. Let
$U'_q(\bb_+)$ be the corresponding Borel subalgebra of $\Uqgp$. The
analogue of Proposition \ref{propJ} holds with $J$ replaced by its
analogue $J'$. The notion of
universal weight function makes sense for the algebra $\Uqgp$ as well.
The conditions (a), (b) remain unchanged, the relation \rf{univ-dec} of
condition (c) is now
\begin{equation}
\begin{split}
&\Delta^{std} \sk{\on{W}_{\bar I}((t_i)_{i\in I})}=
\sum\limits_{I=I_1\sqcup I_2}
\big( \on{W}_{\bar I_i}((t_i)_{i\in I_1})
\otimes
\on{W}_{\bar I_2}((t_i)_{i\in I_2})\big)\times\\&\qquad
\times \big( 1 \otimes \prod_{i\in I_1} \psi^+_{\i(i)}(t_{i})
\big)\times \prod\limits_{\substack{k,l|k<l\\ i_k\in I_1,i_l\in I_2}}
\frac{q^{-(\a_{\i(i_k)},\a_{\i(i_l)})}-t_{i_l}/t_{i_k}}
{1-q^{-(\a_{\i(i_k)},\a_{\i(i_l)})}t_{i_l}/t_{i_k}},
\label{univ-dec1}
\end{split}
\end{equation}
where $\Delta^{std}$ is the coproduct of $U'_q({\mathfrak b}_+)/J'$.

\subsection{Vector-valued weight functions}

Let $V$ be a representation of $\Uqgp$ and $v$ be a vector
in $V$. We call $v$ a {\it singular weight vector} of weight
$\{\l_i(z), i=1,...,r\}$ if
$$
e_{\a}[n]v=0, \quad \psi_{\a_i}^{+}(z)v=\l_i(z)v, \qquad i\in \{1,...,r\},
\quad  n \geq 0,
$$
where $\l_i(z)\in\CC[[z^{-1}]]^\times$. Then $V$ is called a representation
with singular weight vector
$v\in V$ if it is generated by $v$ over $U'_q(\widehat{\mathfrak{g}})$.
It is clear that the ideal $J$ defined by \r{ideal} annihilates any
singular weight vector, i.e., $Jv=0$.

For $V$ a representation generated by the singular weight vector $v$
and $\bar I$ an ordered $\Pi$-multiset, we define a $V$-valued function
\begin{equation}\label{WW77}
w^{\bar I}_V((t_i)_{i\in I})=\on{W}_{\bar I}((t_i)_{i\in I})v
\end{equation}
which is a $V$-valued Laurent formal series. For $I=\emptyset$ we set
$w^{\emptyset}_V=v$. Let $V_i$ ($i=1,2$) be representations generated
by the singular weight vectors $v_i$, with series $\{\l_\a^{(i)}(z) |
\alpha\in \Pi\}$. Then the coproduct property \r{univ-dec1}
of the universal weight function yields the following property of the
$V_1\otimes V_2$-valued function
\begin{equation}
\begin{split}
w^{\bar I}_{V_1\ot V_2}((t_i)_{i\in I})=
\sum\limits_{I=I_1\sqcup I_2}
w^{\bar I_1}_{V_1}((t_i)_{i\in I_1})
\otimes
w^{\bar I_2}_{V_2}((t_i)_{i\in I_2})
\times \\
\times
\prod\limits_{i\in I_1}
\l^{(2)}_{\i(i)}(t_{i})
\ \prod\limits_{\substack{k,l | k<l \\ i_k\in I_1, i_l\in I_2}}
\frac{q^{-(\a_{\i(i_k)},\a_{\i(i_l)})}t_{i_k}-t_{i_l}}
{t_{i_k}-q^{-(\a_{\i(i_k)},\a_{\i(i_l)})}t_{i_l}}\ .
\label{weight1}
\end{split}
\end{equation}
A collection of $V$-valued functions $w_{V}((t_i)_{i\in I})$ for all
possible $\Pi$-multisets $\bar I$, and representations $V$ with
a singular weight vector $v$, is called a vector-valued weight
function or simply a {\it weight function}, if satisfies the relations
\rf{weight1}, the initial condition $w_V(\emptyset)=v$, and depends only on
the isomorphism class of the ordered $\Pi$-multiset $I$.

It is clear how to modify relation \rf{weight1} to define a vector-valued
weight function of the algebra $\Uqg$. Any universal weight function
determines a vector-valued weight function by relation \rf{WW77}.

\subsection{Main Theorems}

Let $\bar I = (I,\prec,\iota)$ be an ordered $\Pi$-multiset.
If $I=\{{i_1},...,{i_n}\}$, with $i_1\prec ... \prec i_n$, we set
\begin{equation} \label{W7}
W_{\bar I}((t_i)_{i\in I}) :=
P\left(f_{\i(i_1)}(t_{i_1})\cdots f_{\i(i_n)}(t_{i_n})
\right).
\end{equation}

\begin{theorem} \label{th1}
The map $\bar I \mapsto W_{\bar I}$ defined by
\rf{W7} is a universal weight function.
\end{theorem}
\noindent
Note that the statement of Theorem \ref{th1} is valid for both algebras
$\Uqg$ and $\Uqgp$.
\medskip

\noindent
{\it Proof}.\ \ First one should check that
$W_{\bar I}((t_i)_{i\in I})\in
U_q({\mathfrak b}_+)((t_{i_n}^{-1}))\ldots ((t_{i_1}^{-1}))$.
This follows from the fact that for any $x\in U_F$
and $\a_i\in\Pi$, there exists an integer $M$, such that
for any $n>M$ we have  $P(xf_{\a_i}[-n])$=0.
We now prove this fact. Define a degree on $U_F$ by
$\on{deg}(h_\alpha[n]) = n$ ($n>0$), $\on{deg}(U_q({\mathfrak h})) = 0$,
$\on{deg}(f_{\alpha_i}[n]) = n$ ($n\in {\mathbb Z}$).
Then $U_F^+,U_f^-\subset U_F$ are graded subalgebras. Moreover,
the nontrivial homogeneous components of $U_F^+$ all have degree $\geq 0$.
We have $U_F = U_F^+ \oplus (U_f^-)^\varepsilon U_F^+$.
Then we take $M$ to be the largest degree of a nonzero homogeneous
component of $x$.

Using the definition of $P$, one also checks that
$$W_{\bar I}((t_i)_{i\in I})
\in t_{i_1}^{-1}U_q({\mathfrak b}_+)((t_{i_n}^{-1}))\ldots ((t_{i_2}^{-1}))
[[t_{i_1}^{-1}]].$$

Let us now show that $\bar I \mapsto W_{\bar I}$ satisfies conditions
(a)--(c) of Section \ref{sectiondef}. Conditions (a) and (b) of Section
\ref{sectiondef} are trivially satisfied. Let us show that (c) is satisfied.
Theorem~\ref{th2} and the following relations
\begin{equation*} 
\begin{split}
\Delta^{(D)} f_\a(z)&= 1\ot f_\a(z)+f_\a(C_2^2z)\ot \psi^+_a(C_2z),
\\
\psi_\a^+(z)f_\b(w)\sk{\psi_\a^+(z)}^{-1}&=
\frac{q^{-(\a,\b)}C^{-1}-w/z}{1-q^{-(\a,\b)}C^{-1}w/z} f_\b(w),
\end{split}
\end{equation*}
imply that relations \rf{univ-dec} are satisfied modulo $\Uqg\ot J$.
But both sides of \rf{univ-dec} belong to $U_q(\bb_+)\ot U_q(\bb_+)$.
Thus they coincide modulo $U_q(\bb_+)\ot J$.
\hfill$\Box$ \medskip

So, we reduced the proof of Theorem 1 to the following statement.
Recall first that $\Delta^{std} \circ P$ defines a map
$U_F \to U_F^+ \otimes U_q({\mathfrak b}_+)$, and $P^{\otimes 2}
\circ \Delta^{(D)}$ defines a map $U_F^+ \to (U_F^+)^{\otimes 2}$.

\begin{theorem}\label{th2} For any element $f\in U_F$
\begin{equation} \label{W8}
\Delta^{std}\big(P(f)\big) \equiv
\left(P\ot P\right)\left(\Delta^{(D)}(f)\right) \quad
 \text{\it{mod}} \quad
 U_F^+ \ot J\,,
\end{equation}
where $J$ is the left ideal of $U_q({\mathfrak b}_+)$ defined by \r{ideal}.
\end{theorem}

\medskip

{\em Proof of Theorem \ref{th2}.} We have $U_F = U_f^- U_F^+$,
therefore\footnote{Recall that if $A\subset U_q(\widehat{\ggg})$,
we set $A^\varepsilon := A \cap \on{Ker}(\varepsilon)$}
$U_F = U_F^+ \oplus (U_f^-)^{\varepsilon} U_F^+$.
So we will prove (\ref{W8}) in the two following cases:
(a) $f\in U_F^+$, and (b) $f = xy$, where $x\in (U_f^-)^\varepsilon$
and $y\in U_F^+$.

Assume first that $f\in U_F^+$. According to (\ref{coid1}),
$\Delta^{(D)}(f)\in U_F \otimes U_F^+$, therefore
$$
(P\otimes P)(\Delta^{(D)}(f)) = (P\otimes \on{id})(\Delta^{(D)}(f)).
$$
According to \rf{twist:rel}, we have
$\Delta^{(D)}(x)={\cal R}_1^{2,1}\Delta^{std}(x)({\cal R}_1^{2,1})^{-1}$
for any $x\in U_q(\widehat{\ggg})$, where $({\cal R}_1^{2,1})^{\pm 1}
\in U_f^-\ot U_e^+$. It follows that
$$
(P\otimes \on{id})(\Delta^{(D)}(f)) = (P\otimes \on{id})
({\cal R}_1^{2,1} \Delta^{std}(f) ({\cal R}_1^{2,1})^{-1}).
$$
Now $\Delta^{std}(f)\in U_F^+ \otimes U_q({\mathfrak b}_+)$, therefore
$\Delta^{std}(f)({\cal R}_1^{2,1})^{-1}\in U_F \otimes U_q({\mathfrak b}_+)$;
it follows
that $(P\otimes \on{id})(\Delta^{std}(f)({\cal R}_1^{2,1})^{-1})$ is
well defined.
Now ${\cal R}_1^{21} \in 1\otimes 1 + (U_f^-)^{\varepsilon} \otimes U_e^+$,
therefore
$$
(P\otimes \on{id})({\cal R}_1^{2,1}\Delta^{std}(f) ({\cal R}_1^{2,1})^{-1})
= (P\otimes \on{id})(\Delta^{std}(f) ({\cal R}_1^{2,1})^{-1}).
$$
Since $({\cal R}_1^{2,1})^{-1} \in 1\otimes 1 + U_f^- \otimes
(U_e^+)^\varepsilon \subset
1\otimes 1 + U_f^- \otimes J$, we have
$$
(P\otimes \on{id})(\Delta^{std}(f) ({\cal R}_1^{2,1})^{-1}) \equiv
(P\otimes \on{id})(\Delta^{std}(f)) \quad \on{mod} \quad U_F \otimes J.
$$
Recall now that $\Delta^{std}(f)\in U_F^+ \otimes U_q({\mathfrak b}_+)$.
This implies that
$$
(P\otimes \on{id})(\Delta^{std}(f)) = \Delta^{std}(f).
$$
Finally, since $f\in U_F^+$, we have $f = P(f)$, therefore
$$
\Delta^{std}(f) = \Delta^{std}(P(f)).
$$
Combining the above equalities and congruence, we get
the desired congruence, when $f\in U_F^+$.

Assume now that $f = xy$, where $x\in (U_f^-)^\varepsilon$ and $y\in U_F^+$.
According to (\ref{coid1}), $\Delta^{(D)}(x)\in U_f^- \otimes U_F$; on the
other hand, the identities $(\varepsilon^{(D)} \otimes \on{id}) \circ
\Delta^{(D)} = \on{id}$ and $\varepsilon = \varepsilon^{(D)}$ imply that
$\Delta^{(D)}(x) \in 1 \otimes x + \on{Ker}(\varepsilon) \otimes
U_q(\widehat{\ggg})$; all this implies that
$\Delta^{(D)}(x)= 1 \otimes x + \sum_i a_i\ot b_i$, where
$a_i\in (U_f^-)^\coun$ and $b_i \in U_F$.

Therefore
$$
(P\otimes \on{id})(\Delta^{(D)}(xy)) =
(P\otimes \on{id}) \big( (1\otimes x + \sum_i a_i \otimes b_i)
\Delta^{(D)}(y)
\big).
$$
Recall that $\Delta^{(D)}(y)\in U_F \otimes U_F^+$. Now since
$a_i\in (U_f^-)^\coun$, we have $P(a_i\eta) = 0$ for any $\eta\in
U_F$, which implies that $(P\otimes \on{id})
((a_i\otimes b_i)\Delta^{(D)}(y)) = 0$. Therefore
$$
(P\otimes \on{id})(\Delta^{(D)}(xy)) =
(P\otimes \on{id}) \big( (1\otimes x) \Delta^{(D)}(y)\big).
$$
Applying $\on{id} \otimes P$ to this identity, we get
$$
(P\otimes P)(\Delta^{(D)}(xy)) =
(P\otimes \on{id}) \circ (\on{id} \otimes P)
\big( (1\otimes x) \Delta^{(D)}(y)\big).
$$
Now since $x\in (U_f^-)^\coun$, we have $P(x\eta) = 0$ for any $\eta\in
U_F$, which implies that $(\on{id} \otimes P)
((1\otimes x)\Delta^{(D)}(y)) = 0$. Therefore
$$
(P\otimes P)(\Delta^{(D)}(xy)) = 0.
$$
Since $P(xy) = 0$, we get
$$
\Delta^{std}(P(xy)) = (P\otimes P)(\Delta^{(D)}(xy)),
$$
which proves the desired congruence for the elements of the form
$xy$, $x\in (U_f^-)^\varepsilon$, $y\in U_F^+$.
\hfill{$\Box$}

\subsection{Functional properties of the universal weight function}

Let $\bar I = (I,\prec,\iota)$ be an ordered $\Pi$-multiset, let $n:= |I|$
and let $\sigma\in{\mathfrak S}_n$ be a permutation. Let $\bar I^\sigma =
(I,\prec_\sigma,\iota)$ be the $\Pi$-multiset such that if
$I = \{i_1,...,i_n\}$, where $i_1\prec ...\prec i_n$, then
$i_{\sigma(1)} \prec_\sigma ... \prec_\sigma i_{\sigma(n)}$.
We call $\bar I^\sigma$ a permutation of $\bar I$.

\begin{proposition} \label{prop:form:W}
There exists a collection of formal functions
$$
\overline W_{n_1,...,n_r}\in (U_q({\mathfrak b}_+)/J)[[u^{(s)}_j |
s \in \{1,...,r\}, j\in \{1,...,n_i\}]],
$$
symmetric in each group of variables $(u^{(s)}_j)_{j = 1,...,n_s}$,
such that
\begin{equation} \label{identity:W}
W_{\bar I}((t_i)_{i\in I}) =
{{\ds \epsilon(\sigma)\prod_{s=1}^r \prod_{i,j\in I_s | i\prec j}
(t_i^{-1} - t_j^{-1})}\over\ds  \prod_{i\in I} t_i
{\prod_{i,j\in I | i\prec j}
(t_i^{-1} - q^{(\iota(i),\iota(j))}t_j^{-1})}}
\overline W_{|I_1|,...,|I_r|}
((t^{-1}_{j_1})_{j_1\in I_1}, ..., (t^{-1}_{j_r})_{j_r\in I_r}),
\end{equation}
where $I_s := \iota^{-1}(\alpha_s)$ (it is an ordered set), and $\sigma
\in {\mathfrak S}_n$ is the shuffle permutation of $\{1,...,n\}$
given by $\{1,...,n\} \simeq I_1\sqcup...\sqcup I_r \to I \simeq
\{1,...,n\}$, where the first and last bijections are ordered,
the order relation
on $I_1\sqcup...\sqcup I_r$ is such that $I_1\prec...\prec I_r$,
and the map $I_1\sqcup...\sqcup I_r \to I$ is such that its restriction to
each $I_s$ is the natural injection.
\end{proposition}

{\em Proof.} For any ordered $\Pi$-multiset $\bar I$, denote by
$F_{\bar I}(t_{i_1},\ldots,t_{i_n})\in U_q({\mathfrak b}_+)
[[t_{i_1}^{\pm 1},...,t_{i_n}^{\pm 1}]]$ the series
$$F_{\bar I}(t_{i_1},\ldots,t_{i_n})=
 f_{\i(i_1)}(t_{i_1})\cdots f_{\i(i_n)}(t_{i_n})\, ,$$
and by $A_{\bar I}(t_{i_1},\ldots,t_{i_n})\in {\mathbb C}[t^{-1}_{i_1},...,
t^{-1}_{i_n}]$ the product
\begin{equation}\label{A}
A_{\bar I}(t_{i_1},\ldots,t_{i_n})=\prod_{1\leq k<l\leq n}(t^{-1}_{i_k}-
q^{({\i(i_k)},{\i(i_l)}}t^{-1}_{i_l})\ .
\end{equation}

Let $\bar I\mapsto G_{\bar I}$ be an assignment taking an ordered
$\Pi$-multiset $\bar I = (I,\prec,\iota)$ (with $I = \{i_1,...,i_n\}$
and $i_1\prec ... \prec i_n$) to $G_{\bar I}(t_{i_1},...,t_{i_n})
\in V[[t_{i_1}^{\pm 1},...,t_{i_n}^{\pm 1}]]$.
We say that $\bar I\mapsto G_{\bar I}$ is antisymmetric, iff
for any $\bar I$ and any $\sigma\in {\mathfrak S}_n$ (where $n = |I|$),
we have $G_{\bar I^\sigma}((t_i)_{i\in I}) = \epsilon(\sigma)
G_{\bar I}((t_i)_{i\in I})$.

The defining relations \rf{1} imply that the assignment
$\bar I \to \overline{F}_{\bar I}$ is antisymmetric, where
$$
\overline{F}_{\bar I}(t_{i_1},\ldots,t_{i_n})=
A_{\bar I}(t_{i_1},\ldots,t_{i_n}) F_{\bar I}(t_{i_1},\ldots,t_{i_n})\,
$$
takes its values in $U_q({\mathfrak b}_+)[[t_{i_1}^{\pm 1},...,
t_{i_n}^{\pm 1}]]$.

Applying $P$, we get that the assignment $\bar I \mapsto \overline W_{\bar I}$
is antisymmetric, where
$$
\overline W_{\bar I}(t_{i_1},...,t_{i_n})
:= A_{\bar I}(t_{i_1},...,t_{i_n}) W_{\bar I}(t_{i_1},...,t_{i_n});
$$
its takes its values in $(U_q({\mathfrak b}_+)/J)[[t_{i_1}^{\pm 1},...,
t_{i_n}^{\pm 1}]]$.

According to the proof of Theorem \ref{th1},
$$
W_{\bar I}(t_{i_1},...,t_{i_n})\in t_{i_1}^{-1}
(U_q({\mathfrak b}_+)/J)((t_{i_n}^{-1}))
...((t_{i_2}^{-1}))[[t_{i_1}^{-1}]];
$$
therefore
$\overline W_{\bar I}(t_{i_1},...,t_{i_n}) :=
A_{\bar I}W_{\bar I}(t_{i_1},...,t_{i_n})$ takes its values in the
same space. The antisymmetry of $\bar I\mapsto \overline W_{\bar I}$
then implies that it takes its values in the intersection of all the
$t_{i_{\sigma(1)}}^{-1}(U_q({\mathfrak b}_+)/J)((t_{i_{\sigma(n)}}^{-1}))
...((t_{i_{\sigma(2)}}^{-1}))[[t_{i_{\sigma(1)}}^{-1}]]$, where $\sigma
\in {\mathfrak S}_n$, i.e., in $$(t_{i_1}...t_{i_n})^{-1}
(U_q({\mathfrak b}_+)/J)[[t_{i_1}^{-1},...,t_{i_n}^{-1}]].$$

If now $V$ is a vector space and $\bar I\mapsto v_{\bar I}(t_{i_1},...,
t_{i_n}) \in V[[t_{i_1}^{-1},...,t_{i_n}^{-1}]]$ is antisymmetric,
one shows that there exists a family of formal series
$\overline v_{n_1,...,n_r}(u^{(1)}_1,...,u^{(r)}_{n_r})\in
V[[u^{(1)}_1,...,u^{(r)}_{n_r}]]$, symmetric in each group of
variables $u^{(i)}_\alpha$ for fixed $i$, such that
$$
v_{\bar I}(t_{i_1},...,t_{i_n}) = \epsilon(\sigma)\prod_{s=1}^r
\prod_{k,l\in I_s | k\prec l} (t_k^{-1} - t_l^{-1})
\overline{v}_{|I_1|,...,|I_r|}((t^{-1}_{i_1})_{i_1\in I_1},...,
(t^{-1}_{i_r})_{i_r\in I_r}).
$$
The result follows. \hfill \qed \medskip

Let $s,t\in \{1,...,r\}$, with $s\neq t$. Let $m:= 1 - a_{\alpha_s\alpha_t}$,
where $(a_{\alpha\beta})_{\alpha,\beta\in\Pi}$ is the Cartan matrix of
${\mathfrak g}$. Let $k_1,...,k_m\in \{1,...,n_s\}$ be distinct, and let
$l\in \{1,...,n_t\}$. Let $H^{st}_{(k_0,...,k_m),l} \subset \oplus_{s=1}^r
{\mathbb C}^{n_s}$ be the subspace of all $(u^{(s)}_j)_{s\in \{1,...,r\},
j\in \{1,...,n_s\}}$, such that
$$
u_l^{(t)} = q^{-{{(m-1)(\alpha_s,\alpha_s)}\over 2}} u^{(s)}_{k_1} =
q^{-{{(m-3)(\alpha_s,\alpha_s)}\over 2}} u^{(s)}_{k_2} = ...
=  q^{{{(m-1)(\alpha_s,\alpha_s)}\over 2}} u^{(s)}_{k_m}.
$$

\begin{proposition} \label{serre:van}
The restriction of $\overline W_{n_1,...,n_r}$ to $H^{st}_{(k_1,...,k_m),l}$
is identically zero.
\end{proposition}

{\em Proof.} The proof is the same as the proof of the similar
statement in \cite{E}, and is based on the quantum Serre relations.
\hfill \qed \medskip

When $q=1$, $W_{\bar I}((t_i)_{i\in I})$ can be computed as follows.
Set
$$
W^{Lie}_{\bar I}((t_i)_{i\in I}) :=
{{[[f_{\iota(i_1)},f_{\iota(i_2)}],...,f_{\iota(i_n)}]^+(t_{i_1})}\over
{(-1 + t_{i_{n-1}}/t_{i_n})...(-1 + t_{i_1}/t_{i_n})}},
$$
where for $x,y\in {\mathfrak g}$, we set $[x,y](t) = [x[0],y(t)]
= [x(t),y[0]]$; then $x(t) = \sum_{n\in{\mathbb Z}} x[n]t^{-n}$, and
$x^+(t) := \sum_{n>0} x[n]t^{-n}$. Then
\begin{equation} \label{classical:W}
W_{\bar I}((t_i)_{i\in I}) =
\sum\limits_{\substack{s\geq 0, (I_1,...,I_s) |
I = I_1\sqcup ...\sqcup I_s, \\ \on{min}(I_1) \prec ... \prec \on{min}(I_s)}}
W^{Lie}_{\bar I_1}((t_i)_{i\in I_1}) ... W^{Lie}_{\bar I_s}((t_i)_{i\in I_s}),
\end{equation}
where the sum is over all the partitions $I = I_1\sqcup ... \sqcup I_s$,
such that $\on{min}(I_1) \prec ... \prec \on{min}(I_s)$, and
$\bar I_i = (I_i,\prec_i,\iota_i)$
is the ordered $\Pi$-multiset induced by $\bar I$.

For $\alpha\in \Pi$, set
where $f_\a^+(z):=\sum_{n>0}\fin{\a}{n}z^{-n}$.

\begin{conjecture} \label{conjA}
$W_{\bar I}((t_i)_{i\in I})$ is a linear combination of
noncommutative polynomials in the $f_\alpha^+(q^k t_i)$, $f_\alpha[0]$
($\alpha\in\Pi$, $k\in{\mathbb Z}$), where the coefficients have the
form $P((t_i)_{i\in I})/\prod_{i,j\in I | i\prec j} (t_i -
q^{-(\iota(i),\iota(j))} t_j)$, and $P((t_i)_{i\in I})$ is a polynomial
of degree $|I|(|I|-1)/2$.
\end{conjecture}

We have $[[f_{\iota(i_1)},f_{\iota(i_2)}],...,f_{\iota(i_n)}]^+(t) =
[[f_{\iota(i_1)}^+(t),f_{\iota(i_2)}[0]],...,f_{\iota(i_n)}[0]]$,
hence the Conjecture is true for any ${\mathfrak g}$ when $q=1$.
According to \cite{KhP1},
it is also true when ${\mathfrak g} = \mathfrak{sl}_2$ or $\mathfrak{sl}_3$
for any $q\neq 0$.
For example, we have in $U_q(\widehat{\mathfrak{sl}}_2)$
(see \cite{KhP1})
$$
\ds{P}\sk{f_\a(t_1)f_{\a}(t_2)}=\ds
f_{\a}^+(t_1)f^+_\a(t_2)-\frac{(q-q^{-1})t_1}{qt_1-q^{-1}t_2}
\sk{f^+_\a(t_1)}^2.
$$
\medskip

{\bf Remark.} ({\ref{classical:W}) is also valid if
$U({\mathfrak n}_-[z,z^{-1}])$
is replaced by $U(\widetilde{\mathfrak n}_-[z,z^{-1}])$, where
$\widetilde{\mathfrak n}_-$ is the free Lie algebra with generators
$f_\alpha,\alpha\in\Pi$; this algebra
is presented by the relations $(z-w)[f_\alpha(z),f_\beta(w)] = 0$
for any $\alpha,\beta\in\Pi$, so the Serre relations do not play a role in the
derivation of (\ref{classical:W}) (where $q=1$). However, the
results of \cite{KhP1} use the quantum Serre relations.
\hfill \qed \medskip

Let now $V$ be a finite dimensional representation of $\Uqg$ with
singular weight vector $v$. Let $\bar I$ be an ordered $\Pi$-multiset
and $w^{\bar I}_V((t_i)_{i\in I})$ be the vector-valued weight function
\begin{equation}
w^{\bar I}_V((t_i)_{i\in I}) :=
P\left(f_{\i(i_1)}(t_{i_1})\cdots f_{\i(i_n)}(t_{i_n})
\right)v\ .
\label{A4}
\end{equation}

\begin{proposition}
Assume that Conjecture \ref{conjA} is true. Then
$w^{\bar I}_V((t_i)_{i\in I})$
is the Laurent expansion of a rational function on ${\mathbb C}^n$.
There exist rational functions
$$
\overline w_{n_1,...,n_r}((u^{(1)}_j)_{
j = 1,...,n_1}, ... , (u^{(r)}_j)_{
j = 1,...,n_r}),
$$
such that the analogue of identity (\ref{identity:W}) holds.

Each function $\overline w_{n_1,...,n_r}$ is symmetric in each group of
variables $(u^{(s)}_j)_{j=1,...,n_s}$. Its only singularities are
poles at $u^{(s)}_j\in S_s$, where $S_s\in{\mathbb C}^\times$
is a finite subset of
${\mathbb C}^\times$. It vanishes on the spaces $H^{st}_{(k_1,...,k_m),l}$.
\end{proposition}

{\em Proof.} According to the theory of Drinfeld polynomials,
the image of $f_\alpha^+(t)$ in $\on{End}(V)$ is a rational
function in $t$ with poles in ${\mathbb C}^\times$. It follows
from Conjecture \ref{conjA} that $w^{\bar I}_V((t_i)_{i\in I})$
is the Laurent expansion of a rational function, which is
regular except for (a) simple poles at $t_i = q^{-(\iota(i),\iota(j))}t_j$,
where $i\prec j$, and (b) poles at $t_i\in S_i$, where $S_i \subset
{\mathbb C}$ is a finite subset.
The form of $W_{\bar I}$ proved in Proposition
\ref{prop:form:W} also implies that $w_V^{\bar I}$ vanishes on the
hyperplanes $t_i = t_j$, where $\iota(i) = \iota(j)$ (as a formal function,
hence as a rational function), and
Proposition \ref{serre:van} implies that
$w_V^{\overline I}$ vanishes on
the spaces $H^{\alpha\beta}_{(k_1,...,k_m),l}$ (as a formal function,
hence as a rational function).

Define then $\overline w_{n_1,...,n_r} := \overline W_{n_1,...,n_r}v$.
Then the analogue of (\ref{identity:W}) holds. It follows from the
properties of $w_V^{\bar I}$ that $\overline w_{n_1,...,n_r}$
is rational, with the announced poles structure.
Since the Laurent expansion of $\overline w_{n_1,...,n_r}$
is symmetric in each group of variables, so is $\overline w_{n_1,...,n_r}$
itself.
\hfill \qed \medskip

\setcounter{equation}{0}
\section{Relation to the off-shell Bethe vectors}

In this section, we relate the universal weight functions to the
off-shell Bethe vectors, in the case of the
quantum affine algebra $\Uqdvaprim$. The algebra $\Uqdvaprim$ is generated
by the modes of the currents $e(z)$, $f(z)$ and $\psi^\pm(z)$.
We will need only the commutation relations between the currents
$f(z),f(w)$ and $f(z),\psi^+(w)$:
\begin{equation}\label{ff}
(qz-q^{-1}w)f(z)f(w)=(q^{-1}z-qw)f(w)f(z)
\end{equation}
\begin{equation}\label{fpsi}
\psi^+(z)f(w)=\frac{q^{-2}-w/z}{1-q^{-2}w/z}f(w)\psi^+(z)
\end{equation}
Using formula \r{ff} we may calculate the projection
$P\sk{f(z_1)\cdots f(z_n)}$.

The algebra $U'_q(\widehat{\mathfrak{sl}}_2)$ also has a realization in terms of
$L$-operators (\cite{RS}):
\begin{equation*}
L^{\pm}(z)
=\sk{\begin{array}{cc}1&
f^\pm(z)\\
0&1\end{array}}\sk{\begin{array}{cc}k^\pm(zq^{-2})^{-1}&0\\
0&k^\pm(z)\end{array}}\sk{\begin{array}{cc}1&0\\
e^\pm(z)&1\end{array}}
=\sk{\begin{array}{cc}A^\pm(z)&B^\pm(z)\\C^\pm(z)&D^\pm(z)
\end{array}}
\end{equation*}
which satisfy
$$
R(u/v)\cdot (L^\epsilon(u)\ot \mathbf{1})\cdot (\mathbf{1}\ot
L^{\epsilon'}(v))=
(\mathbf{1}\ot L^{\epsilon'}(v))\cdot (L^\epsilon(u)\ot
\mathbf{1})\cdot  R(u/v)
$$
with $\epsilon,\epsilon'\in\{+,-\}$, and
\begin{equation*}
\begin{split}
R(z)\ =&\ (qz-q^{-1})\ \left(E_{11}\ot E_{11}+
E_{22}\ot E_{22}\right) +
(z-1)\left(E_{11}\ot E_{11}+E_{22}\ot E_{11}\right)
 +\\& (q-q^{-1})
\left(z E_{12}\ot E_{21}+  E_{21}\ot E_{12}\right)
\end{split}
\end{equation*}
($E_{ij}$ denotes the matrix unit).

According to \cite{DF}, the Gauss coordinates of the $L$-operators
are related to the currents
as follows
$$
e(z)=e^+(z)-e^-(z),\quad f(z)=f^+(z)-f^-(z),\quad \psi^\pm(z)
=\sk{k^\pm(zq^{-2})k^\pm(z)}^{-1}.
$$

Let $v$ be a vector such that $C^+(z)v=0$. The vector-valued function
\begin{equation}\label{BBB}
w(z_1,\ldots,z_n)=B^+(z_1)\cdots B^+(z_n)v
\end{equation}
is called an off-shell Bethe vector. Using the equality
$B^+(z)=f^+(z)k^+(z)$ and the  relation \r{fpsi} we may present the product
\rf{BBB} in terms of the product of the half-currents $f^+(z)$. This gives
the relation
\begin{equation}\label{P-oBv}
B^+(z_1)\cdots B^+(z_n)=\prod_{i<j}^n\frac{qz_i-q^{-1}z_j}{z_i-z_j}
P\sk{f(z_1)\cdots f(z_n)} \prod_{i=1}^n k^+(z_i)
\end{equation}
which shows the relation between the off-shell Bethe vectors and
the weight function \rf{W7}. For a general quantum affine algebra,
the calculation of the  weight functions given by the universal weight
function \rf{W7} is a complicated and interesting problem.
Such calculations for quantum affine algebras
$U_q(\widehat{\mathfrak{sl}}_3)$ and $U_q(\widehat{\mathfrak{sl}}_{N+1})$
are given in \cite{KhP1} and in \cite{KPSLN}.
The relation between the universal weight function \rf{W7} and the
nested Bethe ansatz (\cite{KR83}) will be studied in \cite{KPT}.

\setcounter{section}{6}
\section*{Appendix}
\def\Uqn{{U_{q}^{}({\mathfrak{n}}_+)}}
\def\Uqnn{{U_{q}^{}({\mathfrak{n}}_\pm)}}
\def\Uqb{{U_{q}^{}({\mathfrak{b}}_+)}}
\def\Pim{{\mathrm U}_{{{\mathrm{Im}}}}^+}
\def\tPim{\tilde{\mathrm U}_{{\mathrm{Im}}}^+}
\def\Pimm{{\mathrm U}_{{\mathrm{Im}}}^\pm}
\def\tPimm{\tilde{\mathrm U}_{{\mathrm{Im}}}^\pm}
\def\ne{\tilde{e}}
\def\ng{\tilde{\gamma}}
\def\nga{\tilde{\nu}}
\def\nne{\bar{e}}
\def\nng{\gamma''}
\def\nw{\tilde{w}}
\def\nnw{{w'}}
\def\prt{\ {\prec}\tilde{}\ }
Here we give a proof of the properties \rf{use-pro1} and \rf{use-pro2} of
circular Cartan-Weyl generators, see Section \ref{sectionCW}. The proof uses
the braid group approach to the CW generators, which we describe first.

Let $T_i:\Uqg\to \Uqg$, $i=0,1,\ldots,r$, be the Lusztig automorphisms
\cite{Lus11}, defined by the formulas
\begin{align}\label{Lu1}
T_i(e_{\a_i})&=-e_{-\a_i}k_{\a_i}^{-1}, &
T_i(e_{\a_j})&= \sum\limits_{p+s=-a_{i,j}}(-1)^p q_i^se_{\a_i}^{(p)}
e_{\a_j}e_{\a_i}^{(s)}\,&  
\quad &i\not=j,\\ \label{Lu2}
T_i(e_{-\a_i})&=-k_{\a_i}e_{\a_i}, &
T_i(e_{-\a_j})&= \sum\limits_{p+s=-a_{i,j}}(-1)^p q_i^{-s}e_{-\a_i}^{(s)}
e_{-\a_j}e_{-\a_i}^{(p)}\,&  
\quad &i\not=j
\end{align}
where $e_{\pm\a_i}^{(p)}=e_{\pm\a_i}^{p}/[p]_{q_i}!$.

We attach to the periodic sequence
$\ldots, i_{-1},i_0,i_1,...,i_n,\ldots$ given by \rf{period} the
sequence $(w_n)_{m\in\mathbb{Z}}$ of elements of the Weyl group,
given by  $w_0=w_1=1$, $w_{k+1}=w_{k}s_{i_{k}}$ for $k>0$, and
$w_{l-1}=s_{i_l}w_l$ for $l\leq 0$. Let $\gamma_k$ be the corresponding
positive real roots \rf{gammak}. We have a normal ordering
$\gamma_1\prec\gamma_2\prec...\prec \delta
\prec 2\delta\prec...\prec
\gamma_{-1}\prec\gamma_0$ of the system $\widehat{\Sigma}_+$.

We define real root vectors $e_{\pm\gamma_k}$, where $k>0$ and
$e_{\pm\gamma_l}$, where $l\leq0$ by the relations
\begin{equation}
\label{Ap1}
e_{\pm\gamma_k}=T_{w_k}(e_{\pm\a_{i_k}})\,,
\qquad
e_{\pm\gamma_l}=T_{w_l}^{-1}
(e_{\pm\a_{i_l}})\,,
\end{equation}
that is,
$e_{\pm\gamma_n}=e_{\pm\a_{i_n}}$ for $n=0,1$;\
$e_{\pm\gamma_k}=T_{i_1}T_{i_2}\cdots T_{i_{k-1}}(e_{\pm\a_{i_k}})$ for
$k>1$, and
$
e_{\pm\gamma_l}=T_{i_0}^{-1}T_{i_{-1}}^{-1}\cdots T_{i_{l+1}}^{-1}
(e_{\pm\a_{i_l}})\,$ for $l<0$.
The imaginary root vectors are defined by the relations \rf{im4} and
 \rf{im-rel}. The imaginary root vectors, related to positive roots,
 generate an abelian subalgebra $\Pim\subset\Uqn$. It is characterized by
the properties \cite{Be}
\begin{equation}
\label{Ap8}
p\in\Pim\ \Leftrightarrow\ T_{w_k}^{-1}(p)\in\Uqn
\quad\text{and}\quad
T_{w_l}(p)\in\Uqn\quad\!\text{for all}\quad k>0,\,l\leq 0.
\end{equation}

The root vectors \rf{Ap1}, \rf{im-rel} satisfy the property \rf{use-pro}
(see \cite{Be}) and thus coincide, up to normalization, with the CW generators
of Section \ref{sectionCW}.

Let $c$ be an integer $>0$. Let  $\ldots, j_{-1},j_0,j_1,\ldots$ be the
periodic sequence defined by the rule $j_n=i_{n-c}$
for all $n\in\ZZ$,
$\{\nw_n\}$ the related sequence of elements of the Weyl group, given by
 the rule
$\nw_0=\nw_1=1$, $\nw_{k+1}=\nw_ks_{{j_k}}$ for $k>0$,
and $\nw_{l-1}=s_{{j_l}}\nw_l$ for $l<0$. Let
$\{\ng_n,\ n\in\ZZ\}$ be the corresponding sequence
of real positive roots, $\ng_k=\nw_k(\a_{j_k})$, if $k\geq 0$ and
$\ng_l=\nw_l^{-1}(\a_{j_l})$, if $l\leq 0$.
Let $\{\ne_\gamma\}$ be the CW generators, built by the braid
group procedure, related to the sequence $\{j_k\}$:
$
\ne_{\pm\ng_k}=T_{\nw_k}(e_{\pm\a_{j_k}})$ if $k\geq 1$, and
$\ne_{\pm\ng_l}=T_{\nw_l}^{-1}
(e_{\pm\a_{j_l}})$ if $k\leq 0$. Let ${\tPim}$ be the subalgebra
of $\Uqn$, generated by the imaginary root vectors $\ne_{n\delta}^{(i)}$,
$i=1,...,r$, $n>0$.

 We have the correspondence:
\begin{align}\notag
{\ng_{n}}&=\left\{
\begin{array}{ll}
s_{\a_{i_{1-c}}}\cdots s_{\a_{i_{0}}}({\gamma_{n-c}}),&
 n\not=1,2,...,c,\\
s_{\a_{i_{1-c}}}\cdots s_{\a_{i_{0}}}({-\gamma_{n-c}}),&
n=1,2,...,c.
\end{array}\right.
\\
\label{Ap7}
{\ne}_{\ng_{n}}&=\ T_{{i_{1-c}}}\!\cdots
T_{{i_{0}}}\left(\hat{e}_{\ga_{n-c}}
\right),\qquad n\in\ZZ.
\\
\label{Ap7a}
{\tPim}&=\ T_{{i_{1-c}}}\!\cdots
T_{{i_{0}}}(\Pim)
\end{align}
Indeed, for $n\not= 1,...,c$ we have  ${\ne}_{\ng_{n}}=
T_{{i_{1-c}}}\cdots T_{{i_{0}}}\left({e}_{\ga_{n-c}}\right)$ by the
construction. For $n=1,...,c$ we have ${\ne}_{\ng_{n}}=
T_{{i_{1-c}}}\!\cdots T_{{i_{0}}}\left(T_{w_{n-c}}^{-1}
T_{i_{n-c}}^{-1}T_{w_{n-c}}\right)\big({e}_{\a_{i_{n-c}}}\big)$, which is
equal to
$T_{{i_{1-c}}}\!\cdots T_{{i_{0}}}\left(\hat{e}_{\ga_{n-c}}\right)$
by \rf{Lu1} and \rf{Lu2}. The relation \rf{Ap7a} follows from the
description \rf{Ap8} of the algebra $\Pim$ and its analogue for the algebra
$\tilde{\Pim}$.

Let $\prec$ be the normal ordering of the system $\widehat{\Sigma}_+$, related
to the sequence $\{i_n\}$, $\prec_c$ the corresponding circular order
in $\widehat{\Sigma}$, and
 $\prt$  the normal ordering of $\widehat{\Sigma}_+$, related
to the sequence $\{j_n\}$. For any
 $\tilde{\a},\tilde{\b}\in \widehat{\Sigma}_+$ we have the correspondence:
\begin{equation}\label{Ap9}
\tilde{\a}\prt\tilde{\b}\ \Leftrightarrow\
\a\prec_c\b,
\end{equation}
where ${\a}=s_{\a_{i_0}}\cdots s_{\a_{i_{c-1}}}(\tilde{\a})$, and
${\b}=s_{\a_{i_0}}\cdots s_{\a_{i_{c-1}}}(\tilde{\b})$.

Consider the relation \rf{use-pro} for CW generators, related to
the sequence $\{j_n\}$:
\begin{equation}
\notag
[\ne_{\tilde{\a}},\ne_{\tilde{\b}}]_{q^{-1}}=\sum
 C_{\{n_j\}}^{\{\nga_j\}}(q)\
{\ne}^{n_1}_{\nga_1}{\ne}^{n_2}_{\nga_2}\cdots {\ne}^{n_m}_{\nga_m},
\end{equation}
with $\tilde{\a}\ \prt\ \nga_1\ \prt\ \ldots
\ \prt\
\nga_m\ \prt\ \tilde{\b}$, where
$C_{\{n_j\}}^{\{\nga_j\}}(q)\in {\mathbb C}[q,q^{-1},1/(q^n-1);n\geq 1]$.
 Due to \rf{Ap7}, \rf{Ap7a}, \rf{Ap8},
 automorphism properties of the maps $T_i$, and commutativity of imaginary
root vectors,
this is equivalent to the relation on circular generators:
\begin{equation}
\label{Ap4}
[\hat{e}_{{\a}},\hat{e}_{{\b}}]_{q^{-1}}=\sum
 \bar{C}_{\{n_j\}}^{\{\nu_j\}}(q)\
{\hat{e}}^{n_1}_{\nu_1}{\hat{e}}^{n_2}_{\nu_2}\cdots
{\hat{e}}^{n_m}_{\nu_m},
\end{equation}
with ${\a}\ {\prec}_c\ \nu_1\ {\prec}_c\ \ldots
\ {\prec}_c\ \nu_m\ {\prec}_c\ {\b}$, where
 $\,{\a}=s_{\a_{i_0}}\cdots s_{\a_{i_{c-1}}}(\tilde{\a})$, and
${\b}=s_{\a_{i_0}}\cdots s_{\a_{i_{c-1}}}(\tilde{\b})$;
$\bar{C}_{\{n_j\}}^{\{\nu_j\}}(q)\in {\mathbb C}[q,q^{-1},1/(q^n-1);n\geq 1]$.
This is a particular case of the relation \rf{use-pro1}, when the root
$\a$ satisfies the condition $\,-\delta\ {\prec}_c\ \a\,$ and $\b$ is positive.

Let $d$ be an integer $>0$. Let now $\{j_n\}$ be a periodic
sequence, related to the sequence \rf{period} by the rule $j_n=i_{n+d}$
for all $n$, $\{\nw_n\}$ the related sequence of elements of the Weyl group,
$\{\ng_n,\ n\in\ZZ\}$  the corresponding sequence of real positive roots.
Let $\{\ne_{\pm\gamma}\}$ be   CW generators, built by braid
group procedure, related to the sequence $\{j_k\}$, and
 ${\tPimm}$  the subalgebras
of $\Uqnn$, generated by imaginary root vectors $\ne_{\pm n\delta}^{(i)}$,
$i=1,...,r$, $n>0$.

 We have now the following correspondence:
\begin{align}\label{A10a}
{\ng_{n-d}}&=\left\{
\begin{array}{ll}
s_{\a_{i_{d}}}\cdots s_{\a_{i_{1}}}({\ga_{n}}),&
 n\not=1,2,\ldots d,\\
s_{\a_{i_{d}}}\cdots s_{\a_{i_{1}}}({-\ga_{n}}),&
n=1,2,\ldots d.
\end{array}\right.
\\
\label{Ap10}
{\ne}_{\pm\ng_{n-d}}&=\ T_{{i_{d}}}^{-1}\!\!\cdots
T_{{i_{1}}}^{-1}\left(\nne_{\pm\ga_{n}}
\right),\qquad n\in\ZZ.
\\
\label{Ap11}
{\tPimm}&=\ T_{{i_{d}}}^{-1}\!\!\cdots
T_{{i_{1}}}^{-1}(\Pimm),
\end{align}
where the temporary real root generators $\nne_{-\ng}$ are given by the prescription
$\nne_{\pm\ga_n}=e_{\pm\ga_n}$ for $n\not=1,2,\ldots, d$ and
$\nne_{\pm\ga_n}=T_{w_{n}}T_{i_n}(e_{\pm\a_{i_n}})$ 
 for $n=1,2,\ldots, d$.
Again, the normal ordering ${\prec}\tilde{ }$, attached to
 the sequence $\{j_n\}$, is in accordance with the circular ordering
$\prec_c$:
$
\tilde{\a}\ \prt\ \tilde{\b}\ \Leftrightarrow\
\a\prec_c\b$,
where ${\a}=s_{\a_{i_1}}\cdots s_{\a_{i_{d}}}(\tilde{\a})$, and
${\b}=s_{\a_{i_1}}\cdots s_{\a_{i_{d}}}(\tilde{\b})$.
Consider the relation \rf{use-pro} for CW generators $\ne_{\ng}$, related to
 {\it negative} roots $\tilde{\a}$ and $\tilde{\b}$ :
\begin{equation}
\label{Ap12a}
[\ne_{\tilde{\a}},\ne_{\tilde{\b}}]_{q^{-1}}=\sum
 C_{\{n_j\}}^{\{\nga_{j}\}}(q)\
{\ne}^{n_1}_{\nga_{1}}{\ne}^{n_{2}}_{\nga_{2}}\cdots {\ne}^{m}_{\nga_{m}},
\end{equation}
with $\tilde{\a}, \tilde{\b}, \nga_i, \in-\widehat{\Sigma}_+$,
so that  $-\tilde{\a}\prt-\nga_1\prt\ldots
\prt-\nga_m\prt-\tilde{\b}$.
Due to \rf{Ap10}, \rf{Ap11}, it is equivalent to
\begin{equation}
\label{Ap12}
[\nne_{{\a}},\nne_{{\b}}]_{q^{-1}}=\sum
 \bar{C}_{\{n_j\}}^{\{\nu_j\}}(q)\
{\nne}^{n_1}_{\nu_1}{\nne}^{n_2}_{\nu_2}\cdots {\nne}^{n_m}_{\nu_m},
\end{equation}
with ${\a}\ {\prec}_c\ \nu_1\ {\prec}_c\ \ldots
\ {\prec}_c\
\nu_m\ {\prec}_c\ {\b}$, where  ${\a}=s_{\a_{i_1}}\cdots s_{\a_{i_{d}}}(\tilde{\a})$, and
${\b}=s_{\a_{i_1}}\cdots s_{\a_{i_{d}}}(\tilde{\b})$.
Since all the roots in \rf{Ap12a} are negative, the collection of the
roots $\{\a,\b,\nu_1,\ldots,\nu_m\}$ in \rf{Ap12} contains negative roots
and some  positive roots belonging to the set
$\{\gamma_1,\ldots,\gamma_d\}$. Note, that $\nne_\nu=e_\nu$, if $\nu$ is
negative, and  $\nne_\nu=-k_\nu e_\nu$, if
$\nu\in\{\gamma_1,\ldots,\gamma_d\}$. Now we apply to \rf{Ap12} the
following automorphism of the algebra $\Uqg$:
$$
e_{-\gamma}\mapsto -k_\gamma e_{-\gamma},\qquad
e_{\gamma}\mapsto - e_{\gamma}k_\gamma^{-1},\qquad k_\gamma\mapsto
k_\gamma, \qquad\text{for all}\quad\gamma\in\widehat{\Sigma}_+.$$
One can see, that this automorphism transforms the relation \rf{Ap12}
to the particular case of \rf{use-pro1}:
\begin{equation}
\label{Ap13}
[\hat{e}_{{\a}},\hat{e}_{{\b}}]_{q^{-1}}=\sum
 \tilde{C}_{\{n_j\}}^{\{\nu_j\}}(q)\
{\hat{e}}^{n_1}_{\nu_1}{\hat{e}}^{n_2}_{\nu_2}\cdots
{\hat{e}}^{n_m}_{\nu_m},
\end{equation}
with ${\a}\ {\prec}_c\ \nu_1\ {\prec}_c\ \ldots
\ {\prec}_c\
\nu_m\ {\prec}_c\ {\b}$ and
$\tilde{C}_{\{n_j\}}^{\{\nga_j\}}(q)\in {\mathbb C}[q,q^{-1},1/(q^n-1);n\geq 1]$,
 when the root
$\b$ satisfies the condition $\b\ {\prec}_c\ \delta$ and $\a$ is negative.
The relations \rf{Ap4} and \rf{Ap12} imply \rf{use-pro1} in full
generality. The proof of \rf{use-pro2} is analogous.
\smallskip

{\bf Remark}. There are analogues of all the relations \rf{use-pro},
\rf{use-pro1} and \rf{use-pro2}, in which the order of the products
(equivalently, the order of the root vectors) in the monomials in the right
hand sides are reversed. To derive them, it is sufficient to first
apply the Cartan antiinvolution ${}^*$ to \rf{use-pro}, and then to
apply the arguments of the Appendix to the result.

\section*{Acknowledgements}

This work was supported by the grant INTAS-OPEN-03-51-3350, the
Heisenberg-Landau program, ANR project GIMP N. ANR-05-BLAN-0029-01, the RFBR
grant 04-01-00642 and the RFBR grant for scientific schools NSh-8065.2006.2.
It was partly done during the visits of S.Kh. and S.P. at MPIM (Bonn) in 2004
and 2005, the visit of S.Kh. at CPM (Marseille) in 2006, the visits of S.P. at
LAPTH (Annecy-le-Vieux) in 2006 and at IRMA (Strasbourg) in 2003 and 2006. The
authors wish to thank these centers for their hospitality and stimulating
scientific atmosphere.

\end{document}